\documentclass[a4paper,twoside,11pt]{amsart}

\usepackage[headings]{fullpage}
\usepackage{amsmath} \usepackage{amssymb} \usepackage{amsopn}
\usepackage{amsthm} \usepackage{amsfonts} \usepackage{mathrsfs}
\usepackage{manfnt}
\usepackage[pdfpagelabels,pdftex]{hyperref}
\usepackage[dvips]{graphicx} \usepackage[usenames]{color}
\usepackage[all]{xy}

\hypersetup{
  pdftitle={},
  pdfauthor={},
  pdfsubject={},
  pdfkeywords={},
  colorlinks=true,
  breaklinks=true,
  bookmarksopen=true,
  bookmarksnumbered=true,
  pdfpagemode=UseOutlines,
  plainpages=false}

\DeclareMathOperator{\rA}{A}

\DeclareMathOperator{\rD}{D}

\DeclareMathOperator{\rH}{H}
\DeclareMathOperator{\rI}{I}

\DeclareMathOperator{\rM}{M}

\DeclareMathOperator{\rP}{P}

\DeclareMathOperator{\rR}{R}

\DeclareMathOperator{\rc}{c}

\DeclareMathOperator{\rh}{h}

\newcommand{\bF}{{\mathbb F}}
\newcommand{\bG}{{\mathbb G}}

\newcommand{\bN}{{\mathbb N}}

\newcommand{\bP}{{\mathbb P}}
\newcommand{\bQ}{{\mathbb Q}}
\newcommand{\bR}{{\mathbb R}}

\newcommand{\bZ}{{\mathbb Z}}

\newcommand{\cN}{{\mathscr N}}

\newcommand{\cR}{{\mathscr R}}

\newcommand{\cU}{{\mathscr U}}
\newcommand{\cV}{{\mathscr V}}

\newcommand{\cX}{{\mathscr X}}

\newcommand{\dI}{{\mathcal I}}

\newcommand{\dO}{{\mathcal O}}

\newcommand{\fo}{{\mathfrak o}}

\newcommand{\surj}{\twoheadrightarrow} 
\newcommand{\inj}{\hookrightarrow}

\DeclareMathOperator{\pr}{pr}

\DeclareMathOperator{\Hom}{Hom}

\DeclareMathOperator{\End}{End}

\DeclareMathOperator{\maps}{Maps}

\DeclareMathOperator{\im}{im}

\newcommand{\matzz}[4]{\left(
\begin{array}{cc} #1 & #2 \\ #3 & #4 \end{array} \right)}

\DeclareMathOperator{\Spec}{Spec}

\DeclareMathOperator{\Pic}{Pic}

\newcommand{\redu}{{\rm red}}

\newcommand{\height}{{\rm ht}}

\newcommand{\OO}{\dO}

\newcommand{\Gm}{\bG_m}

\DeclareMathOperator{\pe}{period}
\newcommand{\rig}{{\rm rig}}

\DeclareMathOperator{\Br}{Br}

\DeclareMathOperator{\Frob}{Frob}

\DeclareMathOperator{\Val}{{\rm Val}}

\DeclareMathOperator{\res}{res}

\DeclareMathOperator{\infl}{inf}

\DeclareMathOperator{\Ind}{Ind}

\newcommand{\spez}{{\rm sp}}

\DeclareMathOperator{\Gal}{Gal}

\newcommand{\ep}{\varepsilon}

\newcommand{\topo}{{\rm top}} 
\newcommand{\cons}{{\rm cons}} 
 
\newcommand{\tame}{{\rm tame}} 
\newcommand{\sep}{{\rm sep}}
\newcommand{\alg}{{\rm alg}}
\newcommand{\nr}{{\rm nr}} 
\newcommand{\sh}{{\rm sh}}
\newcommand{\nis}{{\rm Nis}}

\newcommand{\ab}{{\rm ab}} 
\newcommand{\cs}{{\rm cs}} 
 
\newcommand{\tld}{\tilde } 
\newcommand{\vk}{v} 
\newcommand{\fok}{\fo} 
\newcommand{\indlim}{{\rm lim}\kern-16pt\lower5pt\hbox{$\longrightarrow$}}
 
\newcommand{\hhb}[1]{\hbox to#1pt{}}
\newcommand{\kkp}{\kappa}

\newcommand{\tldxs}{\tilde{X}^s}
\newcommand{\valht}{type }
\newcommand{\Valht}{Type }

\newcommand{\bruch}[2]{\genfrac{}{}{0.5pt}{}{#1}{#2}}
\newcommand{\ov}[1]{\mbox{${\overline{#1}}$}} 

\newtheorem{thm}{Theorem}

\newtheorem{prop}[thm]{Proposition}
\newtheorem{lem}[thm]{Lemma}

\newtheorem{cor}[thm]{Corollary}
\newtheorem{conj}[thm]{Conjecture}

\newtheorem*{thm*}{Theorem}
\newtheorem*{thmMAIN}{Main  Result}

\theoremstyle{definition}
\newtheorem{defi}[thm]{Definition}

\theoremstyle{remark}
\newtheorem{rmk}[thm]{Remark}
\newtheorem{nota}[thm]{Notation}

\newtheorem{ques}[thm]{Question}

\newenvironment{pro}[1][Proof]{{\it{#1:}} }{\hfill $\square$}
\newenvironment{pro*}[1][Proof]{{\it{#1:}} }{}
\newenvironment{pro**}[1][]{{\it{#1}} }{\hfill $\square$}

\usepackage{enum}
\newlist{enumer}{enumerate}{2}
\setlist[enumer]{label=(\roman*),align=left,labelindent=0pt,leftmargin=*,widest = (iii)}
\newlist{enumerar}{enumerate}{1}
\setlist[enumerar]{label=\arabic*.,align=left,labelindent=0pt,leftmargin=*,widest = 8.}
\newlist{enumera}{enumerate}{2}
\setlist[enumera]{label=(\arabic*),align=left,labelindent=0pt,leftmargin=*,widest = (8)}
\newlist{enumeral}{enumerate}{2}
\setlist[enumeral]{label=(\alph*),align=left,labelindent=0pt,leftmargin=*,widest = (m)}

\numberwithin{equation}{section}

\begin{document}

\hrule width\hsize
\hrule width\hsize
\hrule width\hsize
\hrule width\hsize

\vspace{1.5cm}

\newcommand{\Ueberschrift}{Arithmetic in the  fundamental group of  a $p$-adic curve \\[1ex] {\small --- On the $p$-adic section conjecture for curves ---} } 

\title[Section conjecture over p-adic fields]{\Ueberschrift} 
\author{Florian Pop}
\address{Florian Pop, Department of Mathematics, University of Pennsylvania, DRL, 209 S 33rd Street, Philadelphia, PA 19104, USA}
\email{pop@math.upenn.edu}
\urladdr{http://www.math.upenn.edu/~pop/}

\author{Jakob Stix}
\address{Jakob Stix, MATCH -- Mathematisches  Institut, Universit\"at Heidelberg, Im Neuenheimer Feld 288, 69120 Heidelberg}
\email{stix@mathi.uni-heidelberg.de}
\urladdr{http://www.mathi.uni-heidelberg.de/~stix/}
\thanks{The first author was supported by the NSF Grant DMS-0801144
and also received support from MATCH of the University of 
Heidelberg. Both authors would like to thank the Isaac Newton Institute 
Cambridge, UK, for the support and the excellent working conditions during 
the NAG Programme in 2009.}

\date{August 11, 2011} 

\maketitle

\begin{quotation} 
  \noindent \small {\bf Abstract} --- 
We establish a valuative version of Grothendieck's section conjecture for curves 
over $p$-adic local fields. The image of every section is contained in the
decomposition subgroup of a valuation which prolongs the $p$-adic valuation
to the function field of the curve.
\end{quotation}


\section{Introduction} \label{sec:intro}

\noindent This note addresses the arithmetic of rational points 
on curves over $p$-adic fields with ramification theory of 
general valuations and the \'etale fundamental group as the 
principal tools.

\subsection{The fundamental group}
Let  $\ov k$ be a 
fixed separable closure of an arbitrary field $k$, and let $\Gal_k:=\Gal(\ov k | k)$ be 
the absolute Galois group of $k$. 

Let $X/k$ be a geometrically connected variety, and  let
$\ov X:=X\times_k\ov k$ be the base change of $X$ to $\ov k$.
The \'etale fundamental group $\pi_1(X, \ov x)$ with base point 
$\ov x$ is an extension
\begin{equation} \label{eq:pi1ext}
1 \to \pi_1(\ov X, \ov{x}) \to \pi_1(X,\ov{x}) \to \Gal_k \to 1,
\end{equation}
where $\pi_1(\ov X,\ov x)$ is the geometric fundamental group 
of $X$ with base point $\ov x$. In the sequel, we  denote the 
extension (\ref{eq:pi1ext}) by 
$\pi_1(X/k)$ and ignore the base points $\ov x$, because 
they will be irrelevant for our discussion.

\subsection{The conjecture}
To a rational point $a \in X(k)$ the functoriality of $\pi_1$ gives rise to a section 
\[
s_a: \Gal_k \to \pi_1(X)
\]
of $\pi_1(X/k)$.  The functor $\pi_1$ depends a priori  on a pointed space  but yields a well defined $\pi_1(\ov X)$-conjugacy class  $[s_a]$ of sections.
The section conjecture of Grothendieck gives a conjectural description 
of the set of all the sections in an arithmetic situation as follows.

\begin{conj}[see Grothendieck \cite{Groth:letter}] 
\label{conj:SC} 
Let $X$ be a smooth, projective and geometrically connected  
curve of genus $\geq 2$ over a number field $k$. Then the map $a \mapsto [s_a]$ 
is a bijection from the set of rational points $X(k)$ onto the 
set of $\pi_1(\ov X)$-conjugacy classes of sections of $\pi_1(X/k)$.
\end{conj} 
Actually, Grothendieck originally made a more general conjecture allowing $k$ to be a finitely generated extension of $\bQ$. Moreover, Grothendieck noticed that $a \mapsto [s_a]$ is injective if $k$ is a 
number field as a consequence of the Mordell-Weil Theorem, 
see \cite{stix:habil} \S10 for details. We refer to \cite{stix:habil} for a general overview on the section conjecture.

\smallskip 

The focus of the present paper is the following local version of Grothendieck's section conjecture.
\begin{conj}[$p$-adic version of the section conjecture] 
\label{conj:padicSC} 
Let $k/\bQ_p$ be a finite extension, and let $X/k$ be a smooth, projective and geometrically connected  
curve  of genus $\geq 2$. Then the map $a \mapsto [s_a]$ 
is a bijection from the set of rational points $X(k)$ onto the set of 
$\pi_1(\ov X)$-conjugacy classes of sections of $\pi_1(X/k)$.
\end{conj} 

To prove the section conjecture 
for curves over a field $k$ as above, it suffices to show that for 
all finite \'etale geometrically connected covers $X'\to X$, if $\pi_1(X'\!/k)$ has a section, then $X'(k)$ is non-empty. This follows by a well known limit argument used already in the work of Neukirch, and introduced to anabelian geometry by Nakamura, while Tamagawa  \cite{tamagawa:affineGC} Prop 2.8 emphasized its significance to the section conjecture, see \cite{Koenigsmann} Lem 1.7 or \cite{stix:periodindex} App.~C.

\subsection{Evidence for the section conjecture}
The first known examples of curves over number fields 
that satisfy the section conjecture were probably given in  
\cite{stix:periodindex} and later \cite{hs}, and also \cite{stix:bm}. 
More recently, Hain \cite{hain:genericcurve} succeeds to verify 
the section conjecture for the generic curve of genus $g \geq 5$. 
These are nevertheless \textit{no sections examples} in the 
sense that there are no sections of $\pi_1(X/k)$ and hence 
no rational points. But as we mentioned above, the 
ostensibly dull case of curves with neither sections nor points is 
exactly the crucial class of examples. 

The  $p$-adic version of the section conjecture has recently 
moved into the focus of several investigations, as pieces 
of evidence for the $p$-adic section conjecture emerged in
recent years. The most convincing piece consists perhaps in 
Koenigsmann's \cite{Koenigsmann} proof of a birational 
analogue of the $p$-adic section conjecture for curves, see 
also  Pop \cite{pop:propsc} for a $\bZ/p\bZ$-metabelian form
of the birational $p$-adic section conjecture. 

On the other
hand, the birational world seems to be quite different from  
the world of curves, because Hoshi \cite{hoshi:propcounterexample} 
showed that a geometrically pro-$p$ version of the section 
conjecture over $p$-adic fields or number fields does not hold, see \cite{stix:habil} \S22.4 for yet more examples of this type.


\subsection{The valuative section conjecture}

Before announcing the main result, 
let us give a valuation theoretic perspective of the 
section conjectures above. 

First, recall that for every 
complete normal curve $X/k$, its closed points $a\in X$ 
are in bijection with the set $\Val_k(K)$ of equivalence 
classes of non-trivial $k$-valuations $w_a$ of the 
function field $K=k(X)$ of $X$ in such a way that the 
residue field $\kappa(w_a)$ of the valuation $w_a$ associated to $a \in X$ equals
$\kappa(a)$. Precisely, the local ring $\OO _{X,a}$ at a closed point $a \in X$
is a discrete $k$-valuation ring of $K$ with 
valuation $w_a\in\Val_k(K)$. Conversely, if $w$ is a 
non-trivial $k$-valuation of $K$, then by the valuative
criterion of properness, the valuation ring $R_w$ of $w$ 
dominates the local ring $\OO_{X,a}$ of a unique point 
$a\in X$. Since $R_w \neq K$, we have $\OO_{X,a} \neq K$ 
and $a\in X$ is a closed point. Hence $R_w=\OO_{X,a}$, 
and $R_w$ is a discrete $k$-valuation ring of $K$. 

Let $\tld K = k(\tld X)$ be the function
field of the universal pro-\'etale cover $\tld X$ of $X$. The extension $\tld K / K$ is Galois with Galois group identified with $\pi_1(X)$. For every $k$-valuation $w$ on $K$ and every prolongation $\tld w$ to $\tld K$ we denote  by $\rD_{\tld w}$ the decomposition subgroup of $\tld w$ in $\pi_1(X)$.
If $w =w_a$ with $a \in X(k)$ then the projection $\rD_{\tld w_a} \to \Gal_k$ is an isomorphism. Hence its inverse gives rise to a section of $\pi_1(X/k)$, the $\pi_1(\ov X)$-conjugacy class of which agrees with $[s_a]$.

\begin{conj}[$\Val_k(K)$ section conjecture] 
\label{conj:valSC}
Let $k$ be a number field or a finite extension of~$\,\bQ_p$.  Let $X/k$ be a smooth, projective and geometrically connected curve $X/k$ of genus $\geq 2$.
Then in the above notations, for every section 
$s:\Gal_k\to\pi_1(X)$ of $\pi_1(X/k)$ there exists a 
valuation $w \in\Val_k(K)$ with $\kappa(w)=k$ 
and a prolongation $\tld w$ of $w$ to $\tld K$ with
$s(\Gal_k)\subseteq \rD_{\tilde w}$.
\end{conj}

Conjecture~\ref{conj:valSC} is in fact equivalent to 
Conjecture~\ref{conj:SC} for $k$ a number field and to 
Conjecture~\ref{conj:padicSC} for $k$ a finite extension of $\bQ_p$.


\subsection{The main result}
Let  $k$ be a finite extension
of $\bQ_p$ with valuation ring $\fok\subset k$ and $p$-adic valuation $\vk$. The 
following richer birational geometric picture unfolds, see 
Appendix~\ref{sec:val} for more details and precise references.

\subsubsection{Geometry}
For a smooth, projective, geometrically 
connected curve $X/k$ we consider the set of all its proper flat
normal models $\cX_i \to \Spec(\fok)$. The set 
$\{\cX_i\}_i$ is partially ordered with respect to the domination 
relation (inducing the identity in $X/k$). We consider 
$ \varprojlim \cX_i$ as an abstract set.
There is a canonical identification 
\[
\Val_\fok(K) = \varprojlim \cX_i
\]
where $\Val_\fok(K)$ is the subspace of the Riemann--Zariski space (see Appendix~\ref{sec:val} for details) of the field $K$ consisting of the $\fok$-valuations, i.e., valuations $w$ whose valuation ring $R_w$ 
satisfies $\fok \subseteq R_w$. Indeed,  for  
$(x_i) \in \varprojlim \cX_i$ the ring 
$R = \indlim\OO_{\cX_i,x_i}$ is a valuation ring of $K$ that contains $\fok$, because the $x_i$ lie on  models over $\Spec(\fok)$.
Conversely, for $w \in \Val_\fok(K)$ with valuation ring $R_w$ the valuative 
criterion of properness yields for every proper model $\cX_i$ a unique point $x_i$ such that $R_w$ dominates $\OO_{\cX_i,x_i}$. The points $x_i$ form a compatible system $(x_i) \in  \varprojlim \cX_i$ and $R_w = \varinjlim \OO_{\cX_i,x_i}$ holds.

\subsubsection{Valuations}
The set $\Val_\fok(K)$ of $\fok$-valuations of $K$ is a disjoint union  
\[
\Val_\fok(K) = \Val_k(K) \amalg \Val_v(K),
\]
where $\Val_v(K)$ is the set of valuations $w$ of $K$ which prolong
$v$ from $k$ to $K$. 
We notice that there is a canonical embedding 
\[
\Val_k(K)\hookrightarrow\Val_v(K)
\]
as follows. Let $w_a\in\Val_k(K)$ be the $k$-valuation, 
corresponding to a closed point $a\in X$. Then the 
residue field $\kappa(w_a)=\kappa(a)$ is a finite extension 
of $k$, hence $\vk$ has a unique prolongation $v_{\kappa(a)}$ 
to $\kappa(w_a)$. The valuation theoretic composition 
$w:=v_{\kappa(a)}\circ w_a$ yields a valuation of $K$ which prolongs $\vk$, thus 
$w$ lies in $\Val_v(K)$. Conversely, if $w\in\Val_v(K)$ is 
a valuation with valuation ring $R$, then $R[1/p]$ is a valuation ring of $K$ that contains $k=\fok[1/p]$. In 
particular, if $R[1/p] \neq K$, then $R[1/p] =  R_{w_a}$ 
for some $w_a \in \Val_k(K)$. For 
$v_{\kappa(a)}$ as above, the valuation 
theoretic composition $v_{\kappa(a)}\circ w_a$ is exactly the valuation $w$ we started 
with. 

We say that $w\in\Val_\fok(K)$ {\bf originates 
from a $k$-rational point}, if there exists $a\in X(k)$ such that
either $w=w_a$ or $w$ is the image of $w_a$ 
under  the canonical embedding $\Val_k(K)\hookrightarrow\Val_v(K)$.

\subsubsection{The main result}
The main result of the present paper is the following positive 
answer to the $\Val_\fok(K)$ variant of the section conjecture
instead of the $\Val_k(K)$ section conjecture above.
\begin{thmMAIN}
Let $k$ be a finite extension of $\hhb2\bQ_p$.
Let $X/k$ be a smooth, projective,  geometrically connected curve of 
genus $\geq 2$, and let $s:\Gal_k \to \pi_1(X)$ be a section of $\pi_1(X/k)$.
\begin{enumera}
\item
There exists an $\fok$-valuation 
$w\in\Val_\fok(K)$ of  $K=k(X)$ and a prolongation $\tld w$ 
to the function field $\tld K=k(\tld X)$ of the universal pro-\'etale cover $\tld X$ such that $s(\Gal_k)$ is contained in the decomposition group $\rD_{\tld w}$ of $\tld w$ in $\pi_1(X)$, see Section~\ref{sec:fixedpt}.
\item
The valuations $w$  to sections $s$ as given by (1) have arithmetic properties as explained in Section~\ref{sec:exotic} and satisfy uniqueness properties as explained in Section~\ref{sec:unique}.
\end{enumera}
\end{thmMAIN}

In Theorem~\ref{thm:sectionsarelocal} we will prove actually a more general assertion concerning  hyperbolic curves $X/k$ that are not necessarily projective.
A smooth, geometrically connected curve $X/k$ 
is called {\bf hyperbolic,\/} if $X$ has negative $\ell$-adic Euler-characteristic 
$\chi(\ov X)=2-2g-r$.
Here $r$ is the number of geometric points needed to smoothly compactify $\ov X$ over $\ov k$ and $g$ is the genus of the smooth compactification. Recall that in characteristic zero, $X$ being hyperbolic is equivalent to $\pi_1(\ov X)$ being non-abelian.

The section conjecture for hyperbolic curves asserts that every conjugacy class of  sections of $\pi_1(X/k)$ is defined as indicated above by a $k$-rational point of the smooth compactification of $X$, or equivalently by a $k$-valuation $v$ of $K$ with residue field $\kappa(v)=k$.

\smallskip

In some sense the Main Result and Theorem~\ref{thm:sectionsarelocal} give an optimal local version of the section conjecture, were Conjecture~\ref{conj:padicSC} to fail. Indeed, if Conjecture~\ref{conj:padicSC} fails, then it fails for a good reason, namely that the projection map $\rD_{\tld w}\to\Gal_k$ from a decomposition subgroup $\rD_{\tld w}\subset\pi_1(X)$ of some valuation $w\in\Val_\fok(K)$ 
admits a section, although the valuation $w$ does not originate from a $k$-rational point. In this respect the Main Result above reduces the $p$-adic section conjecture to a completely local problem, namely to confirm that $\rD_{\tld w}\to\Gal_k$ does not split if $w$ does not originate from a $k$-rational point. The proof of the birational version of the $p$-adic section conjecture as in \cite{pop:propsc} follows the above strategy with $\pi_1(X)$ replaced by $\Gal_K$.

Finally, it was pointed out by Kedlaya that in yet another interpretation of the Main Result above a section of $\pi_1(X/k)$ gives --- if not a $k$-rational point as predicted by Conjecture~\ref{conj:padicSC} --- at least a $k$-\textit{Berkovich} point which is responsible for the section. In light of the above explanations, it remains to be studied, which $k$-Berkovich points might contribute sections of $\pi_1(X/k)$.


\subsection{Relation to a  tempered analogue}

Yves Andr\'e and Shinichi Mochizuki raised the question of relating the main result of the present paper to work of Mochizuki 
\cite{mochizuki:semi} concerning the tempered fundamental group as defined by Andr\'e in \cite{andre:padicavatar}. For a geometrically connected variety $X$ over a finite extension $k$ of $\bQ_p$ the tempered fundamental group  $\pi_1^{\rm temp}(X, \ov x)$ with base point 
$\ov x$ is a pro-discrete group and forms an extension
\begin{equation} \label{eq:temppi1ext}
1 \to \pi_1^{\rm temp}(\ov X, \ov{x}) \to \pi_1^{\rm temp}(X,\ov{x}) \to \Gal_k \to 1
\end{equation}
that we denote by $\pi_1^{\rm temp}(X/k)$, see \cite{andre:padicavatar} \S4. For $X$ a curve, there is a natural inclusion $\pi_1^{\rm temp}(X) \to \pi_1(X)$ that  turns out to be the inclusion into the continuous pro-finite completion of $\pi_1^{\rm temp}(X)$. 
Again, to a $k$-rational point of $X$ we can associate a conjugacy class of sections 
\[
\Gal_k \to \pi_1^{\rm temp}(X).
\] 
The authors learnt from Mochizuki that, based on \cite{mochizuki:semi}, he could prove that a \textbf{tempered section}, i.e., a section $s$ of $\pi_1^{\rm temp}(X/k)$,
always fixes a compatible system $a_i$ of vertices or edges of the dual graphs of stable models $\cX_{i}$ of Galois finite \'etale covers $X_i \to X$.   
The corresponding inductive limit $\tilde{\OO} = \varinjlim_i \OO_{\cX_i,a_i}$ of local rings  $\OO_{\cX,a_i}$ is a local ring of $\tilde{K}$ stabilized by $s(\Gal_k)$. If the answer to our Question~\ref{ques:nonres} is positive, then, in view of Appendix~\ref{sec:center}, the limit $\tilde{\OO}$ is necessarily a valuation ring. 
Nevertheless, at the present state of knowledge, one cannot infer that.
On the other hand, Mochizuki provides an ad hoc argument to indeed find a valuation ring fixed by $s(\Gal_k)$ and dominating $\tilde{\OO}$.

\smallskip

In light of the above, it is natural to ask to what extent \textbf{tempered sections} are different from 
\textbf{pro-finite sections}, i.e., sections of $\pi_1(X/k)$, more precisely: is every section $s: \Gal_k \to \pi_1(X)$ conjugate to a section with image in the subgroup $\pi_1^{\rm temp}(X)$? Although this property of pro-finite sections would easily follow from Conjecture~\ref{conj:padicSC}, we are unable to prove it directly.
In \cite{mochizuki:semi} page 306, Mochizuki speculates that for tempered sections useful arithmetic insights can be provided that are not available in the pro-finite case (in particular the result used in the tempered analogue above). The results of the present note disprove this to some extent.


\subsection{Outline of the paper}
Since the decomposition group $\rD_{\tld w}$ is the stabilizer of the valuation $\tld w$ under the action of $\pi_1(X)$, the property $s(\Gal_k) \subseteq \rD_{\tld w}$ for a section $s$ translates into the existence of a fixed point under the Galois action by $\Gal_k$ via $s$, see Section~\ref{sec:fixedpt}. 

The starting point of our search for a fixed point comes from the Brauer group method, see Section~\ref{sec:brauermethod}, which relies  only on the $\ell$-part of $\pi_1(\ov X)$. The results on the $\ell$-part of $\pi_1(\ov X)$ provided in Sections~\ref{sec:kersp}--\ref{sec:loginertia} 
suffice to show the existence of a fixed point. This is done in Section~\ref{sec:fixedpt} with the help of the  combinatorial Lemma~\ref{lem:combinatorial}. In some sense, the existence of the valuation $\tld w$ fixed by $s(\Gal_k)$  is related to \textit{tame} phenomena.

In Sections~\ref{sec:exotic} and \ref{sec:unique} we address the question about arithmetic and uniqueness properties of the valuations obtained from sections. For the latter we make use 
in a subtle way of the $p$-part in $\pi_1(\ov{X})$,
in particular Tamagawa's non-resolution \cite{tamagawa:resolution}, 
see Section~\ref{sec:bridgesnonres}, 
in order to move apart the $\ell$-parts of inertia groups corresponding to different prime divisors. 
In some sense, uniqueness turns out to be related to \textit{wild} phenomena.

For the convenience of the reader, in Appendix~\ref{sec:val} we provide a complete geometric description of the valuation theory for $p$-adic curves, which is otherwise not sufficiently documented in the literature. Appendix~\ref{sec:hilb} adds a geometric description of the valuation theoretic \textit{Hilbert Zerlegungstheorie}. The notation of the appendix will be used throughout the note.


\setcounter{tocdepth}{1} {\scriptsize \tableofcontents}


\noindent
{\bf Hypothesis.} From now on, if not explicitly stated otherwise, 
we will work under the hypothesis that $k$ is a finite extension 
of $\bQ_p$, hence in particular the residue field $\kappa=\bF$ of the $p$-adic valuation $v$ of $k$ is a finite field. Further, all finite extensions of $k$ are locally compact fields.

\bigskip

\noindent
{\bf Notation.} For notation and terminology of valuation theory and Hilbert Zerlegungstheorie we refer to Appendix~\ref{sec:val} and \ref{sec:hilb}.


\bigskip

\noindent
{\bf Acknowledgements.}
We would like to thank Yves Andr\'e and Shinichi Mochizuki for pointing out potential relations with analoguous questions in the  context of the tempered fundamental group. Especially, in light of Shinichi Mochizuki's comments, our exposition gained in clarity and our pro-finite result was put in the right perspective.
We further like to thank Jordan Ellenberg, H\'el\`ene Esnault, Kiran Kedlaya, Minhyong Kim, Mohamed Sa\"\i di, Tam\'as Szamuely, Akio Tamagawa, and Olivier Wittenberg for their interest in our work.


\section{Detecting inertia of type 1v in the kernel of specialisation}  \label{sec:kersp} 

\subsection{The kernel of  \texorpdfstring{$\spez$}{sp}}
Let $\cX/\fok$ be a model of the smooth, projective, geometrically connected curve $X/k$ as defined in~\ref{sec:model}. The reduced special fibre $Y=\cX_{\bF,\redu}$ is by assumption a strict normal crossing divisor on $\cX$. Let $Y = \bigcup_\alpha Y_\alpha$ be the decomposition into irreducible components $Y_\alpha$ which are smooth, projective curves with field of constants $\bF_\alpha$. The specialisation map of fundamental groups is a surjection 
\[
\spez : \pi_1(X) \surj \pi_1(\cX) = \pi_1(Y)
\]
the kernel of which we denote by $\cN_{X|\cX}$. The inertia group $\rI_w$ of a valuation $w \in \Val_\fok(K)$ lies in $\cN_{X|\cX}$.  Those for valuations of type 1v, see Appendix~\ref{sec:ht1}, generate $\cN_{X|\cX}$ as a pro-finite group by Zariski-Nagata purity of the branch locus.

\subsection{Cohomology on the model}  Let $n \in \bN$ be invertible on $\cX$.
Let $i:Y \inj \cX$ be the closed immersion of the reduced special fibre. Standard computations in \'etale cohomology show $\rR^q i^! \mu_n =0$ for $q=0,1$, and yield the local cycle class map
\begin{equation} \label{eq:0}
\bigoplus_\alpha \underline{\bZ/n\bZ}_{Y_\alpha} \xrightarrow{\sim} \rR^2 i^! \mu_n 
\end{equation}
It follows that $\rH^1_Y\big(\cX,\mu_n\big)$ vanishes and 
\[
\rH^2_Y\big(\cX,\mu_n\big) = \bigoplus_\alpha \rH^0\big(Y_\alpha,\bZ/n\bZ\big).
\]
By proper base change we have $\rH^2\big(\cX,\mu_n\big)=\rH^2\big(Y,\mu_n\big)$, and the relevant part of the localisation sequence reads
\[ 
0 \to \rH^1\big(\cX,\mu_n\big) \to  \rH^1\big(X,\mu_n\big) \xrightarrow{\res} \bigoplus_\alpha \rH^0\big(Y_\alpha,\bZ/n\bZ\big) \xrightarrow{\rho_n} \rH^2\big(Y,\mu_n\big).
\]
Unraveling the definitions for the map $\rho_n$ yields the composite
\[
 \bigoplus_\alpha \bZ/n\bZ \to \Pic(\cX) \otimes \bZ/n\bZ \xrightarrow{\rc_1} \rH^2\big(\cX,\mu_n\big) \xrightarrow{i^\ast} \rH^2\big(Y,\mu_n\big)
\]
which maps $(n_\alpha)$ to $i^\ast \rc_1\big(\OO_\cX(\sum n_\alpha Y_\alpha)\big)$. The map $\res_\alpha$, the $\alpha$ component of $\res$, can be computed by excision and functoriality of the localisation sequence as follows. By abuse of notation, we denote by $\alpha$ also the valuation of type 1v in $\Val_\fok(K)$ corresponding to $Y_\alpha$. Let  $\ov{\alpha}$ be a geometric point  localised in the generic point of $Y_\alpha$. Then, using the notation of 
Appendix~\ref{sec:nearby} and the tame character, see Section~\ref{sec:tch}, we get a commutative diagram with isomorphisms as indicated.
\begin{equation} \label{eq:resa}
\xymatrix@M+1ex@R-1ex@C-1ex{ \rH^1(\pi_1(X),\mu_n) \ar[d] \ar[r]^{\sim} &  \rH^1\big(X,\mu_n\big) \ar[d] \ar[r]^{\res} &  \rH^2_Y\big(\cX,\mu_n\big) \ar[d] \ar@{=}[r] & \bigoplus_\alpha \rH^0\big(Y_\alpha,\bZ/n\bZ\big) \ar[d]^{\pr_\alpha} \\
 \rH^1\big(\rI_\alpha,\mu_n\big) \ar[r]^{\sim} _{\infl} & \rH^1\big(\cU_\alpha^\sh,\mu_n\big) \ar[r]^{\sim} & \rH^2_{{\ov{\alpha}}}\big(\cX^\sh_\alpha,\mu_n\big) \ar@{=}[r] & \bZ/n\bZ}
\end{equation}
The inflation map $\infl$ in the diagram is an isomorphism because the map $\pi_1(\cU_\alpha^\sh) \surj \rI_\alpha$ is an isomorphism on the prime to $p$ part due to enough tame ramification along each $Y_\alpha$, see Proposition~\ref{prop:Ilog} (3) below. Consequently, the map $\res_\alpha$ is essentially the map induced by restriction from $\pi_1(X)$ to $\rI_\alpha$. 

\subsection{Cohomology of the special fibre} 
Let $\ov{Y} = Y \times_\bF \bF^\alg$ be the geometric reduced special fibre. Let $\dI_{\alpha,n}$ be the permutation module $\bZ/n\bZ[\Hom_\bF(\bF_\alpha,\bF^\alg)]$ as a $\Gal_\bF= \Gal(\bF^\alg/\bF)$-module. 
The degree maps of the components describe a $\Gal_\bF$-equivariant  isomorphism
\[
\rH^2\big(\ov{Y},\mu_n\big) = \bigoplus_\alpha \dI_{\alpha,n}.
\]
The relevant cohomology of $Y$ computes via the Leray spectral sequence as
\[
0 \to \rH^1\big(\bF,\rH^1(\ov{Y},\mu_n)\big) \to \rH^2\big(Y,\mu_n\big) \xrightarrow{(\deg_\alpha)} \bigoplus_\alpha \rH^0\big(\bF, \dI_{\alpha,n}\big) \to 0,
\]
where $\deg_\alpha$ is the degree map on the component $Y_\alpha$.  If we fix an $\bF$-embedding $\bF_\alpha \subset \bF^\alg$, then 
\[
\dI_{\alpha,n} = \Ind_{\bF_\alpha}^{\bF}\big(\bZ/n\bZ\big) \]
becomes canonically isomorphic to the induced module with respect to  $\Gal_{\bF_\alpha} \subset \Gal_\bF$.

\subsection{Unramified extensions of the base}

Now we perform the limit of the above computations over unramified extensions $k'/k$ viewing the result as $\Gal_\bF = \pi_1(\fok)$-modules. In other words, we take the stalk at $\Spec(\bF^\alg) \to \Spec(\fok)$ of the higher direct images for $\cX/\fok$. The unramified base changes do not destroy the good properties that $\cX$ has by assumption as a model, 
see Appendix~\ref{sec:model}, and no modification by blow-ups is necessary. We get an exact sequence of $\Gal_\bF$-modules as follows.
\begin{equation} \label{eq:1}
0 \to \rH^1\big(\ov{Y},\mu_n\big) \to  \rH^1\big(X \times_k k^\nr,\mu_n\big) \to 
 \bigoplus_\alpha \Ind_{\bF_\alpha}^{\bF} \! \big(\bZ/n\bZ\big) \xrightarrow{\ov{\rho}_n} \bigoplus_\alpha \Ind_{\bF_\alpha}^{\bF} \! \big(\bZ/n\bZ\big).
\end{equation}
The map $\ov{\rho}_n = (\deg_\alpha) \circ \rho_n$ is a matrix with entries from $\End(\bZ/n\bZ) = \bZ/n\bZ$ with rows and columns indexed by the $\Gal_\bF$-set of irreducible components of $\ov{Y}$. This is nothing but the intersection matrix for the reduced geometric special fibre modulo $n$. 

\subsection{\texorpdfstring{$\ell$}{ell}-adic coefficients}

The local cycle class (\ref{eq:0})  is compatible with change of coefficients $\mu_n \subset \mu_{nd}$ with $p \nmid d$ via the commutative diagram
\[
\xymatrix@M+1ex@R-1ex{ \bigoplus_\alpha \underline{\bZ/n\bZ}_{Y_\alpha}  \ar[d]^{\cdot d} \ar[r]^\cong & \rR^2 i^! \mu_n \ar[d] \\
\bigoplus_\alpha \underline{\bZ/nd\bZ}_{Y_\alpha}  \ar[r]^\cong & \rR^2 i^! \mu_{nd}.
}\]
Taking the direct limit of (\ref{eq:1}) for $n=\ell^r$, $r \geq 0$ we obtain 
an exact sequence of $\Gal_\bF$-modules
\begin{equation} \label{eq:2}
0 \to \rH^1\big(\ov{Y},\bQ_\ell/\bZ_\ell(1)\big) \to  \rH^1\!\big(X_{k^\nr},\bQ_\ell/\bZ_\ell(1)\big) \to 
 \bigoplus_\alpha \Ind_{\bF_\alpha}^{\bF} \! \big(\bQ_\ell/\bZ_\ell\big) \xrightarrow{\ov{\rho}} \bigoplus_\alpha \Ind_{\bF_\alpha}^{\bF} \! \big(\bQ_\ell/\bZ_\ell\big)
\end{equation}
Here $\ov{\rho}$ is a matrix with entries from $\End(\bQ_\ell/\bZ_\ell) = \bZ_\ell$ with  rows and columns indexed by the $\Gal_\bF$-set of irreducible components of $\ov{Y}$, which is the intersection matrix for the reduced geometric special fibre, and moreover takes values in $\bZ \subseteq \bZ_\ell$. 
As an integral matrix, the matrix of  $\ov{\rho}$ is symmetric, negative semi-definite with radical given by the rational mutiples of the divisor of the special fibre with its multiplicities, see Mumford \cite{mumford:surface} \S1.

\subsection{Unramified extensions of the model}

We compute the limit of (\ref{eq:2}) over all finite \'etale covers $\cX'$ of $\cX$. The comments on the preservation of the good properties of the model still hold true, so we can use   (\ref{eq:2}) for all covers.  With $X'$ the generic fibre and $Y'$ the special fibre of $\cX'$ as above we have 
\[
\varinjlim_{\cX'}  \rH^1\big(\ov{Y'},\bQ_\ell/\bZ_\ell(1)\big) = 0 
\]
and
\[
\varinjlim_{\cX'} \rH^1\big(X' \times_k k^\nr,\bQ_\ell/\bZ_\ell(1)\big) = \rH^1\big(\cN_{X|\cX},\bQ_\ell/\bZ_\ell(1)\big)
\]
by compatibility of cohomology of pro-finite groups and discrete coefficients with limits. 
If $\cX'$ corresponds to an open subgroup $H \subset \pi_1(\cX)$, then 
the part of $\rH^2_{Y'}\big(\cX',\bQ_\ell/\bZ_\ell(1)\big)$ due to components of $Y'$ above $Y_\alpha$ is given by 
\[
\maps_{\pi_1(Y_\alpha)}\big(\pi_1(Y)/H,\bQ_\ell/\bZ_\ell\big) =  \left[\maps_{\pi_1(Y_\alpha)}\big(\pi_1(Y),\bQ_\ell/\bZ_\ell\big) \right]^H.
\]
In the limit over all $\cX' \to \cX$ we obtain the smooth induction 
\[
\Ind_{\pi_1(Y_\alpha)}^{\pi_1(Y)}\big(\bQ_\ell/\bZ_\ell\big) = \bigcup_H   \left[\maps_{\pi_1(Y_\alpha)}\big(\pi_1(Y),\bQ_\ell/\bZ_\ell\big) \right]^H.\]
The transfer maps in the limit $\varinjlim \rH^2\big(\ov{Y}',\bQ_\ell/\bZ_\ell(1)\big)$ multiply by the respective degrees. All  components of positive genus are dominated by components with degree an arbitrary high power of $\ell$, even abelian covers, as $\pi_1^\ab (\ov{Y}_\alpha) \otimes \bZ_\ell  \inj \pi^\ab_1 (\ov{Y} \otimes \bZ_\ell)$ shows. Therefore in the limit  only the components $Y_\beta$ of genus $g_\beta=0$ survive.  Alltogether, we get the sequence
\begin{equation} \label{eq:3}
0 \to   \rH^1\big(\cN_{X|\cX},\bQ_\ell/\bZ_\ell(1)\big)  \to 
 \bigoplus_\alpha \Ind_{\pi_1(Y_\alpha)}^{\pi_1(Y)}\big(\bQ_\ell/\bZ_\ell\big) 
  \xrightarrow{\cR}  \bigoplus_{\beta, \ g_\beta = 0}  \Ind_{\pi_1(Y_\beta)}^{\pi_1(Y)}\big(\bQ_\ell/\bZ_\ell\big).
\end{equation}
Let $\rI^{\ab,\ell}_\alpha = \rI^\ab_\alpha \otimes \bZ_\ell$ be the $\ell$-Sylow group of $\rI_\alpha^\ab$.
Taking Pontrjagin duality with Tate-twist, i.e, $\Hom\big(-,\bQ_\ell/\bZ_\ell(1)\big)$,  and using (\ref{eq:resa}) we get the exact sequence of pro-finite $\pi_1(Y)$-modules
\begin{equation} \label{eq:4}
\bigoplus_{\beta, \ g_\beta = 0} \rI_\beta^{\ab,\ell}[[{\pi_1(Y)}/{\pi_1(Y_\beta)}]] \xrightarrow{\cR^\vee(1)} \bigoplus_\alpha  \rI_\alpha^{\ab,\ell} [[{\pi_1(Y)}/{\pi_1(Y_\alpha)}]] \to  \cN_{X|\cX}^\ab \otimes \bZ_\ell \to 0.
\end{equation}
Here we have used the notation $M[[G/G_0]]$ for a finitely generated $\bZ_\ell$-module $M$ and a closed subgroup $G_0$ of a profinite group $G$ to denote
\[
M[[G/G_0]] = \varprojlim_{H} M \otimes_{\bZ_\ell} \bZ_\ell[H\backslash G/G_0],
\]
where $H$ ranges over the open normal subgroups of $G$ and $ \bZ_\ell[H\backslash G/G_0]$ is the permutation module on the set $H\backslash G/G_0$ with coefficients in $\bZ_\ell$. The dual of the induced module $\Ind_{G_0}^G(\bQ_\ell/\bZ_\ell)$ equals $\bZ_\ell[[G/G_0]]$ due to the identification
\[
(\bZ/\ell^n\bZ)[H\backslash G/G_0] = \Hom\big(\maps_{G_0}(G/H,\bruch{1}{\ell^n}\bZ/\bZ),\bQ_\ell/\bZ_\ell\big)
\]
mapping $HgG_0 \in H\backslash G/G_0$ to the evaluation $f \mapsto f(g^{-1})$ for $f \in \maps_{G_0}(G/H,\bruch{1}{\ell^n}\bZ/\bZ)$. 

The composition of $\cR^\vee(1)$ in (\ref{eq:4}) with the projection to the part of components of genus zero yields a map, which is a projective limit indexed over finite \'etale covers $Y' \to Y$ of maps as follows
\[
\bigoplus_{\beta, \ g_\beta = 0}  \rI_\beta^{\ab,\ell}[[{\pi_1(Y')} \backslash {\pi_1(Y)}/{\pi_1(Y_\beta)}]] \to
\bigoplus_{\beta, \ g_\beta = 0} \rI_\beta^{\ab,\ell}[[{\pi_1(Y')} \backslash {\pi_1(Y)}/{\pi_1(Y_\beta)}]] .
\]
For each $Y' \to Y$ the map is given by a matrix with the intersection pairing of $Y' \subset \cX'$ restricted to the genus zero components. If in $Y$ at least one component has genus at least $1$, then this matrix is negative definite,  see Mumford \cite{mumford:surface} \S1, and hence  the map is injective and remains so in the projective limit over all $Y'$.

\begin{prop} \label{prop:detect1v}
Let  $\alpha_1,\ldots,\alpha_r \in \Val_\fok(K)$ be valuations of type 1v which belong to distinct components of $Y$ with positive genus. Then the natural map 
\[
\bigoplus_{i=1}^r \rI_{\alpha_i}^{\ab,\ell} \inj  \cN_{X|\cX}^\ab \otimes \bZ_\ell \]
is injective and  the $\ell$-Sylow subgroups of any two distinct  $\rI_{\alpha_i}$ intersect  trivially in $\pi_1(X)$.
\end{prop}
\begin{pro}
The computation above shows that the images of the natural map
\[
\bigoplus_{i=1}^r \rI_{\alpha_i}^{\ab,\ell}  \to 
\bigoplus_\alpha \rI_\alpha^{\ab,\ell}[[{\pi_1(Y)}/{\pi_1(Y_\alpha)}]] \]
and of $\cR^\vee(1)$ meet only trivially. 
\end{pro}


\section{The logarithmic point of view towards inertia} \label{sec:loginertia}


\subsection{The tame character}  \label{sec:tch}
Let $\Gamma_w$ be the value group of $w \in \Val_v(K)$. Let $\hat{\bZ}'(1)$ be the prime to $p$ Tate module of roots of unity in the separable closure  $\kappa(w)^\sep$ of $\kappa(w)$. The \textbf{tame character} at $w$ is the surjective homomorphism
\[
\chi: \Gal_{K_w^\sh} \surj \Hom\big(\Gamma_w,\hat{\bZ}'(1) \big)
\]
that maps $\sigma \in \Gal_{K_w^\sh}$ to the homomorphism 
\[\chi_\sigma: \gamma \mapsto \left(\sigma(\sqrt[n]{t_\gamma})/\sqrt[n]{t_\gamma}\right)_n\]
with $t_\gamma$ being an arbitrary element of value $w(t_\gamma) = \gamma$.
The kernel of the tame character $\chi$ is the $p$-Sylow group of $\rI_w = \Gal_{K_w^\sh}$.


\subsection{Enters the logarithmic fundamental group} A model $\cX$, see Appendix~\ref{sec:model}, can be naturally equipped with a log-regular fs-log structure by the divisor $\cX_{\bF,\redu}$.  We obtain a quotient
\[
\pi_1(X,\bar{\eta}) \surj \pi_1^{\log}(\cX)
\]
which has a tractable group structure by a logarithmic van Kampen theorem applied to the logarithmic special fibre. 

\begin{prop} \label{prop:logblowup}
(1) For a map $f:\cX' \to \cX$ between models of $X$ the induced map 
\[
\pi_1^{\log}(f) \, : \ \pi_1^{\log}(\cX') \to \pi_1^{\log}(\cX)
\]
is an isomorphism. 

(2) Let $X$ admit a stable model $\cX_{\rm stable}$. Then $\cX_{\rm stable}$ admits an fs log structure which is log regular,
and for any model $\cX$ of $X$ the natural map $f: \cX \to \cX_{\rm stable}$  the induced map 
\[
\pi_1^{\log}(f) \; : \ \pi_1^{\log}(\cX) \to \pi_1^{\log}(\cX_{\rm stable})
\]
is an isomorphism. 
\end{prop}
\begin{pro}
In both cases $f$ is a composition of blow-ups which can be enriched to logarithmic blow-up maps. A logarithmic blow-up map yields an isomorphism of log fundamental groups by \cite{fujiwarakato:logetale} 2.4, see also \cite{illusie:overview} Thm 6.10 or \cite{stix:diss} Cor 3.3.11.
\end{pro}


\subsection{Logarithmic inertia groups}

We denote by $\rI_w^{\log}$ (resp.\ $\rI_y^{\log}$) the  image of $\rI_w$ (resp.\ $\rI_y$) in $\pi_1^{\log}(\cX)$, which is a pro-finite group of order prime to $p$.
The log structure on $\cX$ induces a log structure on $\cX_y^{\sh}$ for every  $y$. The group 
$\rI_y^{\log}$ is nothing but the image of 
\[
\pi_1^{\log}(\cX_y^\sh,\bar{\xi}_y) \to \pi_1^{\log}(\cX,\bar{\eta}) .
\]

\begin{lem} \label{lem:tamecoord}
Let $y$ be a geometric point of $\cX$ which lies above a point of the special fibre. 
\begin{enumera}
\item
The natural map 
\[
\pi_1^{\log}(\cX_y^\sh,\bar{\xi}_y) \to  \Hom\Big(\OO^\ast(\cU_y^\sh)/\OO^\ast(\cX^\sh_y),\hat{\bZ}'(1) \Big)
\]
induced by the tame character is an isomorphism.
\item \label{enum:2}
Let $y$ lie over the generic point of the component $Y_\alpha$ of the special fibre associated to a valuation $\alpha$ of type 1v. Then we have canonically 
\[
\hat{\bZ}'(1) =   \Hom\Big(\OO^\ast(\cU_\alpha^\sh)/\OO^\ast(\cX^\sh_\alpha),\hat{\bZ}'(1) \Big) = \pi_1^{\log}(\cX_\alpha^\sh,\bar{\xi}_\alpha) 
\]
\item Let $y$ lie over a closed point of the special fibre.  Then the canonical map 
\[
\bigoplus_{\alpha \; ; \ y \in Y_\alpha} \hat{\bZ}'(1) \to \pi_1^{\log}(\cX_y^\sh,\bar{\xi}_y),
\]
given essentially by $\ref{enum:2}$ above and  restriction of units is an isomorphism.
Here $\alpha$ ranges over the valuations associated to components $Y_\alpha$ of the special fibre with $y \in Y_\alpha$. 
\end{enumera}
\end{lem}
\begin{pro}
This standard result in log geometry follows from Abhyankar's Lemma and Zariski--Nagata purity of the branch locus, see \cite{illusie:overview} Example 4.7 or \cite{stix:diss} Cor 3.1.11.
\end{pro}

The following lemma describes the behaviour of logarithmic inertia groups under changes of the model.

\begin{lem} \label{lem:blowup}
Let $f:\cX' \to \cX$ be a blow-up of a closed point of the model $\cX$, above which we have the geometric point $y$. Let $\alpha$ denote the valuation associated to a component $Y_\alpha$ of the special fibre of $\cX$ which contains $y$, and let $\ep$ denote the valuation associated to the exceptional divisor $E$ of the blow-up on $\cX'$. Let $z$ (resp.\ $x$) be geometric points of $\cX'$ which lift $y$ and lie on $E$, such that $z$ also lies on the strict transform of $Y_\alpha$ (resp.\ such that $x$ lies in the smooth locus of the reduced special fibre of $\cX'$). The map $f$ yields maps 
\[
f_{z,y} : \pi_1^{\log}(\cX_z^{',\sh},\bar{\xi}_z) \to \pi_1^{\log}(\cX_y^\sh,\bar{\xi}_y)  \quad \text{ and } \quad  f_{x,y} : \pi_1^{\log}(\cX_x^{',\sh},\bar{\xi}_x) \to \pi_1^{\log}(\cX_y^\sh,\bar{\xi}_y)
\]
with the following description in terms of the canonical coordinates provided by Lemma~\ref{lem:tamecoord} \ref{enum:2}.
\begin{enumera}
\item Let $y$ be a node of the reduced special fibre of $\cX$ with the other component through $y$ besides $Y_\alpha$ being the component $Y_\beta$ associated to the valuation $\beta$.
Then $f_{z,y}$ is the isomorphism 
\[
\matzz{1}{1}{0}{1} : \hat{\bZ}'(1) \oplus \hat{\bZ}'(1) \xrightarrow{\sim}  \hat{\bZ}'(1) \oplus \hat{\bZ}'(1) 
\]
with respect to the ordering $(\alpha,\ep)$ and $(\alpha,\beta)$. And $f_{x,z}$ is the diagonal injection 
\[
\hat{\bZ}'(1) \inj \hat{\bZ}'(1) \oplus \hat{\bZ}'(1).
\]
\item  Let $y$ be a smooth point of the reduced special fibre of $\cX$. Then $f_{z,y}$ is the surjection given by the sum
\[
 \hat{\bZ}'(1) \oplus \hat{\bZ}'(1) \surj \hat{\bZ}'(1), 
\]
and $f_{x,y}$ is the identity isomorphism 
\[
\hat{\bZ}'(1) \xrightarrow{\sim} \hat{\bZ}'(1).
\]
\end{enumera}
\end{lem}
\begin{pro}
It all comes down to compute the valuations of the pull back to $\cX'$ of local parameters at $y$. Let the center of the blow up be the ideal $(u,v)$ with $u=0$ describing $Y_\alpha$ and, if present $v=0$ describing $Y_\beta$. Then near $z$ we find that $v=0$ describes $E$ while $u/v=0$ describes the stric transform of $Y_\alpha$. Hence $\ep(u) = \ep(v) = 1$ which leads to the matrix in (1). The remaining calculations are of the same kind but simpler.
\end{pro}

Let $\rI_k$ (resp.\ $\rI_k^{\tame} = \hat{\bZ}'(1)$) be the inertia (resp.\ tame inertia) group of $\Gal_k$ (resp.\ its tame quotient $\Gal_k^{\tame}$).  The projection $\pi_1(X,\bar{\eta}) \to \Gal_k$ maps the inertia (resp.\ the log inertia) groups associated to points or valuations to $\rI_k$ (resp.\ $\rI_k^\tame$).

\begin{prop} \label{prop:Ilog}
Let $w \in \Val_\fok(K)$ be a valuation. The structure of $\rI_w^{\log}$ is as follows.
\begin{enumera}
\item If $w$ is of type 1h, then  $\rI_w^{\log}=1$.
\item If $w$ is of type 2h with $w=v_y \circ \alpha$, then $\rI_w^{\log} = \hat{\bZ}'(1)$ with the natural map $\rI_w^{\log} \to \rI_k^\tame$ being multiplication by the ramification index $e(v_y/v)$.
\item If $w$ is of type 1v, then $\rI_w^{\log} = \hat{\bZ}'(1)$ with the natural map $\rI_w^{\log} \to \rI_k^\tame$ being multiplication by the ramification index $e(w/v)$.
\item If $w$ is a valuation of \valht $2$ and on some model $\cX$ the center $x_w$ lies only on one component of the special fibre associated to the valuation $\alpha$, then 
$\rI_w^{\log} \subseteq \rI_\alpha^{\log} = \hat{\bZ}'(1)$. This applies in particular to valuations of type 2h, 2u$_{\rm sm}$ and 2u$_{\rm alt}$.
\item If the center of a valuation $w$ of \valht $2$ is a node of the reduced special fibre on all models, then 
\[
\hat{\bZ}'(1) \oplus \hat{\bZ}'(1) \surj \rI_w^{\log} = \langle  \rI_\alpha^{\log}, \rI_\beta^{\log} \rangle
\]
where $\alpha,\beta$ are the valuations of type $1v$ which correspond to the components through the center of $w$ on a given model $\cX$ of $X/k$, and approriate base points have been chosen.
\end{enumera}
\end{prop}
\begin{pro}
All assertions follow from Proposition~\ref{prop:logblowup} and Lemma~\ref{lem:blowup}.
\end{pro}

\smallskip

We would like to stress, that if in (5) in fact we have equality $\hat{\bZ}'(1) \oplus \hat{\bZ}'(1) = \rI_w^{\log}$, then the decomposition as a direct sum depends on the choice of model according to Lemma~\ref{lem:blowup} (1). 

\begin{prop} \label{prop:inertiadesing}
Let $X/k$ have stable reduction $\cX_{\rm stable}/\fok$. Let $f:\cX \to \cX_{\rm stable}$ be the minimal regular resolution of the stable model. Let $y \in \cX_{\rm stable}$ be a node with singularity of type $A_n$ in the intersection of the two distinct components $Y_{\alpha_1},Y_{\alpha_2}$, such that $f^{-1}(y)$ equals a chain  of irreducible divisors $E_1, \ldots,E_{n-1}$, that links the strict transforms $E_0,E_n$ of $Y_{\alpha_1},Y_{\alpha_2}$, i.e., such that 
\begin{enumer}
\item $E_i$ meets $E_{i-1}$ and $E_{i+1}$ each in a single node for $i = 1,\ldots,n-1$,
\item and $E_i \cong \bP^1$ for $i = 1,\ldots,n-1$.
\end{enumer}
Let $\ep_i$ be the valuation of type 1v associated to $E_i$ for $i = 1,\ldots,n-1$.  Then for $i = 1,\ldots,n-1$ the natural map 
\[
\bZ_\ell(1) = \rI_{\ep_i}^{\log} \otimes \bZ_\ell \to \rI_y^{\log} \otimes \bZ_\ell  \subset \Big(\rI_{\alpha_1}^{\log} \oplus \rI_{\alpha_2}^{\log} \Big) \otimes \bQ_\ell = \bQ_\ell(1) \oplus \bQ_\ell(1)
\]
is given by multiplication with $(\bruch{i}{n},\bruch{n-i}{n})$ and thus remains injective after projection to each component.
\end{prop}
\begin{pro*}
\'Etale locally around $y$ the situation is as follows. The local ring is $R=\fok[u,v]/(uv-\pi^n)$ with $u$ and $v$ parameters along $Y_{\alpha_1}$ and $Y_{\alpha_2}$ and $n$ is the thickness of the double point singularity, that equals the length of the chain of $E_i$'s connecting the strict transforms $E_0,E_n$ in the minimal resolution $f: \cX \to \cX_{\rm stable}$. 

The map $f$ enhances to a log blow up and thus has a combinatorial description within the fs monoid $Q = \ov{\rM}_y = R[1/\pi]^\ast/R^\ast$ which is spanned by $u,v$ and $\pi$. We give a description of the dual monoids because in the end $\Hom(Q,\hat{\bZ}'(1))$ equals the tame inertia $\rI_y^{\log}$ at $y$. In coordinates dual to $u,v$ we find 
\[
Q^\vee = \{(a,b) \in \bruch{1}{n}(\bN_0)^2 \ ; \  a+b \in \bZ \}.
\]
The log blow up corresponds to a subdivision of $Q^\vee$ as follows.  The component $E_i$ comes from the dual of the submonoid $P_i ^\vee \subset Q^\vee$ generated by $(\bruch{i}{n},\bruch{n-i}{n})$ and hence the node  $E_i \cap E_{i+1}$ is given by the dual of 
\[
\langle (\bruch{i}{n},\bruch{n-i}{n}), (\bruch{i+1}{n},\bruch{n-i-1}{n}) \rangle \subset Q^\vee,
\]
for $i=0,\ldots, n-1$.
From the fact, that this monoid is isomorphic to $\bN^2$ we see again that $\cX$ is indeed regular in the nodes $E_i \cap E_{i+1}$.  Moreover, using the special values $i = 0$ and $i = n-1$ 
it follows that indeed
\[
\rI_y^{\log} \otimes \bZ_\ell = Q^\vee \otimes \bZ_\ell(1) \subset  \Big(\rI_{\alpha_1}^{\log} \oplus \rI_{\alpha_2}^{\log} \Big) \otimes \bQ_\ell = \bQ_\ell(1) \oplus \bQ_\ell(1)
\]
with respect to the identity map on coordinates 
$(a,b) \in Q^\vee \otimes \bZ_\ell(1)  \mapsto (a,b) \in \bQ_\ell(1) \oplus \bQ_\ell(1)$.
The proposition now follows from the identification 
\[
\ \hspace{2cm} \rI_{\ep_i}^{\log} \otimes \bZ_\ell = P_i^\vee \otimes \bZ_\ell(1) = (\bruch{i}{n},\bruch{n-i}{n}) \cdot \bZ_\ell(1) \subset Q^\vee \otimes \bZ_\ell(1) =  \rI_y^{\log} \otimes \bZ_\ell. \hspace{2cm} \qed
\]
\end{pro*}


\begin{cor} \label{cor:inertiawrtstablemodel}
Let $X$ admit a stable model $\cX_{\rm stable}$, and let $\cX$ be a model of $X$ with natural map $f: \cX \to \cX_{\rm stable}$. Let $\tilde{w}$ be a prolongation of the valuation $w \in \Val_\fok(K)$ of \valht $2$, and let $y \in \cX_{\rm stable}$ be the image $f(x_w)$ of the center of $w$ under $f$. 
The logarithmic inertia $\rI_{\tilde{w}|w}^{\log}$ is a subgroup of $\rI_y^{\log}$ and the intersection 
\[
\rI_{\tilde{w}|w}^{\log} \cap \bigoplus_{\alpha} \rI_{\tilde{\alpha}|\alpha}^{\log}
\]
is of finite index in $\rI_{\tilde{w}|w}^{\log}$, where $\alpha$ ranges over the valuations of type 1v associated to irreducible components of the special fibre of $\cX_{\rm stable}$ that contain $y$,  
and the $\tilde{\alpha}$ are prolongations to $\tilde{K}$ such that $\rI_{\tilde{\alpha}|\alpha}^{\log} \subset \rI_y^{\log}$. More precisely, on every log \'etale cover of $\cX_{\rm stable}$ the center of $\tilde{\alpha}$ is determined as the generic point of a component of the special fibre passing through the point above $y$ which determines the embedding $\rI^{\log}_{y} \subset \pi_1^{\log}(\cX)$.
\end{cor}
\begin{pro}
This all follows from Proposition~\ref{prop:logblowup} and Proposition~\ref{prop:inertiadesing} except for the claim that $ \bigoplus_{\alpha} \rI_{\tilde{\alpha}|\alpha}^{\log}$ with $\tilde{\alpha}|\alpha$ as in the statement is indeed a subgroup in $\rI_y^{\log}$. The latter is a consequence of Proposition~\ref{prop:detect1v} and Lemma~\ref{lem:visible} below.
\end{pro}


\subsection{Visible valuations of type 1v}

We first recall a useful property of the dual graph of the special fibre in the tower of all finite \'etale covers.


\subsubsection{Disentangling the dual graph} \label{sec:disentangle}

For a proper model $\cX/\fok$ the reduced special fibre $Y = \cX_{\bF,\redu}$ has a dual graph $\Gamma=\Gamma_Y$ which describes the combinatorics of the components $Y_\alpha$ of $Y$ with their mutual intersection. The completely split fundamental group $\pi_1^\cs(Y)$ is the quotient of $\pi_1(Y)$ which describes mock covers of $Y$, i.e., those finite \'etale covers which are geometrically completely split over every generic point of $Y$. 

\begin{lem} \label{lem:disentangle}
Two different components of the cover of $Y$ corresponding to the maximal geometrically abelian exponent $2$ quotient of $\pi_1^\cs(Y)$ intersect at most once.
\end{lem}
\begin{pro}
This is simply topological covering theory of finite graphs.
\end{pro}

\begin{rmk}
Lemma~\ref{lem:disentangle} says that any two components of the reduction \textbf{disentangle} after taking a  finite \'etale cover, even with good completely split reduction, which by definition means that any chosen preimages intersect at most once.
\end{rmk}


\subsubsection{Visible valuations} \label{sec:visible}

\begin{defi}
A valuation of type 1v is \textbf{visible} if the associated component of the special fibre is covered by a component of positive genus  in the reduction of a model of a suitable finite \'etale cover of the generic fibre. In this case we also call the associated component visible. A valuation of type 1v is \textbf{invisible} if it is not visible.
\end{defi}

\begin{ques} \label{ques:nonres}
An old question of the first author that resists all our efforts to resolve it asks whether all components are visible. 
\end{ques}

\begin{rmk}
(1) 
The important work of Tamagawa  \cite{tamagawa:resolution} towards this Question~\ref{ques:nonres} only guarantees that  over any closed point $y \in \cX_{\rm stable}$  of the stable model any fine enough model $f: \cX \to \cX_{\rm stable}$ has non-stable rational components in 
$f^{-1}(y)$ which are visible, see also Section~\ref{sec:bridgesnonres}.

(2) For a given  finite \'etale cover of the generic fibre there are only finitely many components of the special fibre which have positive genus. As $\pi_1(\ov{X},\bar{\eta})$ is topologically finitely generated, the category of all finite \'etale covers has up to scalar extension only countably many isomorphism classes. Hence, there are at most countably many visible components on each model. As in our case the residue field of $k$ is finite, every $X/k$ admits only countably many components of the special fibres of its models up to taking strict transforms, so that we have no cardinality reason to argue that not all components are visible. 
\end{rmk}

Let $\alpha$ be a valuation of type 1v and let $Y_\alpha$ be its associated component of the reduced special fibre $\cX_{\bF,\redu}$ of a model $\cX$ of $X/k$.  When we endow $\cX$ with the vertical log structure coming from the special fibre $Y \inj \cX$ and moreover $Y_\alpha$ by the induced log structure, then we know from Lemma~\ref{lem:disentangle}  and \cite{Mz:globalprofinite} Prop 4.2 or  \cite{stix:diss} Prop 6.2.11, that 
\[
\pi_1^{\log}(Y_\alpha) \inj \pi_1^{\log}(\cX,\bar{\eta}) 
\]
is injective, whenever $Y_\alpha$ is a \textbf{stable} component, i.e, the strict transform of a component from  the stable model of $X$.
In particular every stable component acquires positive genus in a finite Kummer \'etale logarithmic cover of $\cX$. The lemma below is an immediate  consequence. 

\begin{lem} \label{lem:visible}
For a valuation $\alpha$ of type 1v the following are equivalent.
\begin{enumeral}
\item $\alpha$ is visible.
\item $Y_\alpha$ is dominated by a component of the stable model of a finite \'etale cover of the generic fibre.
\item The genus of a component from a model of a finite \'etale cover of the generic fibre which dominate $Y_\alpha$ tends to infinity in a cofinal tower of all finite \'etale covers of the generic fibre.
\hfill $\square$
\end{enumeral} 
\end{lem}


\subsection{The kernel of (log-)specialisation prime-to-\texorpdfstring{$p$}{p}}

\begin{prop} \label{prop:compare}
Let $j:X \subset \cX$ be the inclusion of a smooth, projective curve $X/k$ as an open subset a model $\cX$. The natural map 
\[
\cN_{X|\cX} \to \cN^{\log}_{X|\cX}
\]
induced by $\pi_1(j): \pi_1(X) \to \pi_1^{\log}(\cX)$ from the kernel $\cN_{X|\cX} = \ker\big(\spez: \pi_1(X) \surj \pi_1(\cX)\big)$ to the kernel 
$\cN^{\log}_{X|\cX}$ of the \em{`forget log map'} $\pi_1(\ep): \pi_1^{\log}(\cX) \surj \pi_1(\cX)$ 
induces an isomorphism on the maximal prime-to-$p$ quotients.
\end{prop}
\begin{pro}
The log scheme $\cX$ is log regular and has $X$ as its locus of trivial log structure. The map $\pi_1(j): \pi_1(X) \to \pi_1^{\log}(\cX)$ induces an isomorphism on maximal prime-to-$p$ quotients by Fujiwara--Kato's purity for the log fundamental group \cite{fujiwarakato:logetale} Thm 3.1, see also \cite{illusie:overview} Thm 7.6. 
The proposition follows in the limit from this reasoning applied to arbitrary finite \'etale covers $\cX' \to \cX$ and the corresponding generic fibre $X' \subset \cX'$.
\end{pro}


\section{Sections and the Brauer group method} \label{sec:brauermethod}


\subsection{Local--global principles for the Brauer group} 
Lichtenbaum constructs for a smooth projective curve $X$ over the $p$-adic field $k$ a perfect  duality pairing \cite{lichtenbaum:duality} \S5 Thm 4 
\begin{equation} \label{eq:lipairing}
\Br(X) \times \Pic(X) \to \Br(k) = \bQ/\bZ.
\end{equation}
The vanishing of the left kernel of  (\ref{eq:lipairing}) translates into the injectivity of the map
\begin{equation} \label{eq:lgpoints}
\Br(X) \inj \prod_{a \in X_0} \Br(\kappa(a)),
\end{equation}
which evaluates a Brauer class on $X$ in every closed point $a \in X$. 

Let $\cX$ be a model of $X/k$. A closed point $a \in X$ has the henselian local scheme $Z_a \subset \cX$ as its Zariski-closure in $\cX$. The unique closed point of $Z_a$ is the topological intersection $y_a = Z_a \cap \cX_\bF$ with the special fibre of the model. Because $Z_a$ is henselian, its inclusion to $\cX$ lifts to the scheme of nearby points $\cX^{\rh}_{y_a}$, and the lift induces a map $\Spec(\kappa(a)) \to \cU^{\rh}_{y_a}$ that lifts the point $a \in X$. 
It follows immediately from (\ref{eq:lgpoints}) that we also have a local global principle 
\begin{equation} \label{eq:lgmodel}
\Br(X) \inj \prod_{y \in \cX_{0}} \Br(\cU_y^{\rh})
\end{equation}
with respect to the nearby points $\cU_y^{\rh}$ associated to the closed points $y\in \cX_0$ of a model.
In the direct limit over all models of $X/k$ we find the composition
\begin{equation} \label{eq:lgallmodels}
\Br(X) \inj \varinjlim_{\cX} \prod_{y \in \cX_{0}} \Br(\cU_y^{\rh}) \to \prod_{w \in \Val_v(K)} \Br(K_w^{\rh}).
\end{equation}
The last map follows  from the restriction map
\begin{equation} \label{eq:brlimit}
  \varinjlim_{\cX} \Br(\cU^{\rh}_{x_w}) \xrightarrow{\sim}  \Br(U_w^{\rh}) \inj \Br(K_w^{\rh}),
\end{equation}
and is injective by purity for the Brauer group. Moreover, if $w$ is not  of type 2h,  the map (\ref{eq:brlimit}) is an isomorphism by the compatibility of henselisation and the Brauer group with direct  limits.

\begin{prop} \label{prop:locus}
Let $\cX$ be a model of $X/k$. Let $A \in \Br(X)$ be a Brauer class.
Then the set 
\[
\{y \in \cX_{\bF,\redu} \  ; \ A \text{ is nontrivial in } \Br(\cU_y^{\rh})\}
\]
is closed in the constructible topology $\cX_{\cons}$.
\end{prop}
\begin{pro}
We only have to argue that if $A$ vanishes in $\Br(\cU_\alpha^{\rh})$ for some generic point $\alpha$ of $\cX_{\bF,\redu}$, then $A$ vanishes in $\Br(\cU_y^{\rh})$ for all but finitely many closed points $y$ in the closure of $\alpha$. But if $A$ vanishes in $\Br(\cU_\alpha^{\rh})$, then this occurs already on $V = \cV \times_\cX X$ for some strict \'etale neighbourhood $\cV \to \cX$ of $\alpha$. For allmost all $y$ in question  the natural map $\cU_y^{\rh} \to X$ factors over $V$ and therefore $A$ also vanishes in those $ \Br(\cU_y^{\rh})$.
\end{pro}

\begin{cor} \label{cor:locus}
Let $A \in \Br(X)$ be a Brauer class. Then the set 
\[
\{w \in \Val_\fok(K) \  ;  \ A \text{ is nontrivial in } \Br(K_w^{\rh})\}
\]
is closed in the patch  topology on $\Val_\fok(K)$.
\end{cor}

\begin{pro}
Corollary~\ref{cor:locus} follows at once from Proposition~\ref{prop:locus} because $\Val_\fok(K)$ with the patch topology is homeomorphic to  $\varprojlim_{\cX} \cX_{\cons}$.
\end{pro}

Proposition~\ref{prop:locus} and the fact that a projective limit of nonempty compact spaces is nonempty shows that the composite map in (\ref{eq:lgallmodels}) is also injective 
\cite{pop:diss} Thm 4.5. 
More precisely, by taking limits and exploiting  the compactness of the patch topology we find the following generalization of (\ref{eq:lgallmodels}):
\begin{thm}[\cite{pop:diss} Thm 4.5] \label{thm:popdiss}
Let $M/k$ be a function field of transcendence degree $1$ over $k$. Then  the following restriction map is injective
\begin{equation} \label{eq:lgp}
\Br(M) \inj \prod_{w|v} \Br(M_w^{\rh}),
\end{equation}
where the product ranges over all valuations $w$ of $M$ extending the $p$-adic valuation $v$ on $k$. 
\end{thm}


\subsection{Computation of the Brauer group of a henselisation --- case of \valht 2}  \label{sec:compbr}
We are interested in controlling the kernel of $\Br(k) \to \Br(K_w^{\rh})$ for a valuation $w \in \Val_v(K)$. We will compute for each model $\cX$ a relevant subgroup of $\Br(\cU_y^{\rh})$ and take the limit  over all models as in (\ref{eq:brlimit}).

Let $y \in \cX$ be a closed point of a model. The cohomology sheaves with support of $\Gm$ for $i: Y_y^{\rh} := \cX^{\rh}_y \setminus \cU_y^{\rh} \inj \cX^{\rh}_y$ are $\rR^q i^! \Gm = 0$ for $q=0,2$  and 
\[
\rR^1i^! \Gm = j_\ast\Gm/\Gm \cong  \bigoplus i_{Y_\alpha,\ast} {\bZ}
\]
with the isomorphism induced by the valuation $w_\alpha$ on the function field $K$ of $\cX$ defined by the components $i_\alpha : Y_\alpha \inj \cX_y^{\rh} $ of $Y^{\rh}_y$. Moreover, $i_\ast \rR^3 i^!\Gm =  \rR^2 j_\ast \Gm$ with open immersion $j:\cU_y^{\rh} \subset \cX_y^{\rh}$ has stalk $(\rR^3 i^!\Gm)_{\bar{y}} = \Br(\cU_y^{\sh})$ and therefore the map
\[
\Br(\cU_y^{\rh}) \to \rH^3_{Y_y^{\rh}}\big(\cX_y^{\rh},\Gm\big) \to \rH^0\big(Y_y^{\rh},\rR^3 i^!\Gm\big) = \rH^0\big(y,(\rR^3 i^!\Gm)_{\bar{y}}\big)
\]
has kernel the relative Brauer group $\Br(\cU_y^{\sh}/\cU_y^{\rh})$ of classes in $\Br(\cU_y^{\rh})$ which die when restricted to $\cU_y^{\sh}$. In the limit over all models we get $\Br(U_w^{\sh}/U_w^{\rh}) \subset \Br(K_w^{\rh})$.
By the computation of $\Br(k)$ along $i: \Spec(\bF) \inj \Spec(\fok)$ as 
\[
\Br(k) = \rH^3_{\Spec(\bF)}\big(\Spec(\fok),\Gm\big) = \rH^2\big(\bF,\rR^1i^! \Gm\big) = \Hom\big(\Gal_\bF,v(k) \otimes \bQ/\bZ\big) \cong \bQ/\bZ
\]
the subgroup $\Br(\cU_y^{\sh}/\cU_y^{\rh})$ receives the image of the restriction map $\Br(k) \to \Br(\cU_y^{\rh})$.
Let $\Gal_y$ be the absolute Galois group of the residue field $\kappa(y)$ at $y$.
Using $\Br(\cX_y^{\rh}) = \Br(\kappa(y)) = 0$ and $\rH^3\big(\cX_y^h,\Gm\big) = \rH^3\big(\kappa(y),\Gm\big) = 0$, the relative cohomology sequence for $(\cX_y^{\rh},\cU_y^{\rh})$ yields an isomorphism of $\Br(\cU_y^{\sh}/\cU_y^{\rh})$ with 
\[
\rH^2\big(Y_y^{\rh},\rR^1i^! \Gm\big) =  \bigoplus_\alpha \Hom\big(\Gal_y, w_\alpha(K) \otimes \bQ/\bZ\big) = \Hom\big(\Gal_y, \Gm(\cU^{\rh}_y)/\Gm(\cX^{\rh}_y)\big).
\]
In the limit over all models we obtain the following proposition:
\begin{prop}
Let $w \in \Val_v(K)$ be a valuation of \valht $2$. The map $\Br(k) \to \Br(U_w^{\sh}/U_w^{\rh})$ is isomorphic to the map:
\begin{enumera}
\item  If $w$ is not of type 2h, then
\[
\Hom\big(\Gal_\bF,v(k) \otimes \bQ/\bZ\big) \to \Hom\big(\Gal_{\kappa(w)},w(K) \otimes \bQ/\bZ\big),
\]
\item If $w = v_{\kappa(a)} \circ w_a$ is of type 2h refining $w_a$ of type 1h for a closed point $a \in X$ and $v_{\kappa(a)}$ the $p$-adic valuation of the residue field $\kappa(a)$, then
\[
\Hom\big(\Gal_\bF,v(k) \otimes \bQ/\bZ\big) \to  \Hom\big(\Gal_{\kappa(w)},v(\kappa(a)) \otimes \bQ/\bZ\big),
\]
\end{enumera}
where the maps are defined by the inclusion map $v(k) \inj w(K)$, resp.\ $v(k) \inj v(\kappa(a))$, of value groups and the restriction map $\Gal_{\kappa(w)} \to \Gal_\bF$. 
\hfill $\square$
\end{prop}

\begin{cor} \label{cor:criterion}
The class of invariant $1/\ell$ survives in $\Br(K_w^{\rh})$ for a valuation $w$ of \valht $2$ if and only if the degree of the residue field extension $\kappa(w)/\bF$ is prime to $\ell$ and the value $w(\pi)$ of a uniformizer $\pi$ of $k$ is not divisible by $\ell$ in the value group $w(K)$. \hfill $\square$
\end{cor}


\subsection{Local--semilocal principle for the Brauer group}

Let $\alpha$ be a valuation of type 1v, and let $Y_\alpha$ be the associated divisor in suitably fine models. The scheme $Y_\alpha$ is a smooth projective curve over a finite extension $\bF_\alpha/\bF$ as field of constants. We define $\Br'(K_\alpha^{\rh})$ as the preimage of
 $\rH^1(\pi_1(Y_\alpha), \bQ/\bZ) \subset  \rH^1(\kappa(\alpha),\bQ/\bZ) = \rH^2(\kappa(\alpha),\rR^1 i^!\Gm)$ under the natural map
\[
\Br(K_\alpha^{\rh}) \to \rH^3_{\alpha}\big(\cX_\alpha^{\rh},\Gm\big),
\]
from the relative cohomology sequence and the local to global spectral sequence associated to the regular embedding $i: \Spec(\kappa(\alpha)) \inj \cX_\alpha^{\rh}$. The subgroup $\Br'(K_\alpha^{\rh})$ receives the image of the restriction map $\Br(k) \to \Br(K_\alpha^{\rh})$.
By $\rH^3\big(\cX_\alpha^{\rh},\Gm\big) = \rH^3\big(\kappa(\alpha),\Gm\big) = 0$, we can extract from the relative cohomology sequence the exact sequence 
\begin{equation} \label{eq:valsequence}
0 \to \Br(\kappa(\alpha)) \to \Br'(K_\alpha^{\rh}) \to \rH^1\big(\pi_1(Y_\alpha),\bQ/\bZ\big) \to 0.
\end{equation}
A valuation $w_{y\alpha} = v_y \circ \alpha$ which refines $\alpha$ by means of a closed point $y \in Y_\alpha$ has henselisation $K_{w_{y\alpha}}^{\rh}$ containing $K_\alpha^{\rh}$. The restriction map $\Br(K_\alpha^{\rh}) \to \Br(K_{w_{y\alpha}}^{\rh})$ respects the respective subgroups $\Br'(K_{w_{y\alpha}}^{\rh}) \to \Br(U_{w_{y\alpha}}^{\sh}/U_{w_{y\alpha}}^{\rh})$. The value group of $w$ sits in an exact sequence of torsion free groups
\[
0 \to v_y(\kappa(\alpha)) \to w_{y\alpha}(K) \to w_\alpha(K) \to 0,
\]
which therefore remains exact after applying $\Hom(\Gal_y, (-)\otimes \bQ/\bZ)$. 
The restriction maps on Brauer groups for all such $w_{y\alpha} := v_y \circ w_\alpha$ fit into a map
\begin{equation} \label{eq:semilgp}
\xymatrix@R-2ex@C-3ex{ 0 \ar[r] & \Br(\kappa(\alpha)) \ar[r] \ar[d]^(0.4){\lambda_1} & \Br'(K_\alpha^{\rh}) \ar[r] \ar[d]^(0.4){\lambda_2} & \Hom(\pi_1(Y_\alpha),\bQ/\bZ) \ar[r] \ar[d]^(0.4){\lambda_3} & 0\\
\lower -8pt\hbox{$0$} \ar@<0.8ex>[r] & {\displaystyle  \prod_{y \in Y_\alpha}}  \Hom\big(\Gal_y,  \bQ/\bZ\big) \ar@<0.8ex>[r] &    { \displaystyle \prod_{y \in Y_\alpha}   \Br(U_{w_{y\alpha}}^{\sh}/U_{w_{y\alpha}}^{\rh})} \ar@<0.8ex>[r] & {\displaystyle \prod_{y \in Y_\alpha} \Hom\big(\! \Gal_y, \bQ/\bZ\big)  } \ar@<0.8ex>[r] & \lower -8pt\hbox{$0$}}
\end{equation}
 of exact sequences.
The homomorphism $\lambda_1$ is injective by the local--global principle for the Brauer group of the function field $\kappa(\alpha)$. The homomorphism $\lambda_3$ restricts an unramified character to the decomposition group of $y \in Y_\alpha$ and is injective because the set of Frobenius elements is dense in $\pi_1(Y_\alpha)$. Hence by the $5$-Lemma we deduce the following local--semilocal principle:
\begin{prop} \label{prop:ht2even}
The restriction map
\begin{equation} \label{eq:lslp}
\Br'(K_{\alpha}^{\rh}) \inj  \prod_{y \in Y_\alpha}   \Br(U_{w_{y\alpha}}^{\sh}/U_{w_{y\alpha}}^{\rh})
\end{equation}
is injective. \hfill $\square$
\end{prop}


\subsection{The Brauer group of the decomposition pro-cover of a section}

In this section we fix a section $s:\Gal_k \to \pi_1(X,\bar{\eta})$ of the fundamental group extension $\pi_1(X/k)$.


\subsubsection{The decomposition pro-cover of a section} 
The section $s$ induces a right action of $\Gal_k$ on the universal pro-\'etale cover $\tilde{X}$ of the curve $X/k$. The corresponding quotient 
\[
\tldxs = \tilde{X}/s(\Gal_k)
\]
 is the maximal subcover $X'/X$ of $\tilde{X}/X$ such that the section $s$ lifts to a section of the composition $\pi_1(X') \subset  \pi_1(X,\bar{\eta}) \to \Gal_k$. In fact, $\pi_1(\tldxs)$ is nothing but the image $s(\Gal_k)$ in $ \pi_1(X,\bar{\eta})$.


\subsubsection{Relative Brauer groups and sections} The \textbf{relative Brauer group} $\Br(X/k)$  of the $p$-adic curve $X/k$ is the kernel of the pullback map $\Br(k) \to \Br(X)$. By a theorem of Roquette and Lichtenbaum, see \cite{lichtenbaum:duality} Thm p.120, we know that $\Br(X/k)$ is cyclic of order the index of $X$. The index of $X$ is defined as $\gcd(\deg(D))$, where $D$ ranges over all $k$-rational divisors on $X$.

\begin{prop} \label{prop:brauerxs}
For a section $s$ of $\pi_1(X/k)$ and any $\ell \not= p$ the natural map 
\[ 
\Br(k) \otimes \bZ_\ell \xrightarrow{\sim} \Br(\tldxs) \otimes \bZ_\ell 
\]
is an isomorphism.
\end{prop}
\begin{pro}
The Leray spectral sequence yields an exact sequence
\begin{equation} \label{eq:LerayGm}
0 \to \Br(X/k) \to \Br(k) \to \Br(X) \to \rH^1\big(k,\Pic_X\big).
\end{equation}
For every finite subcovers $X' \to X$ of $\tldxs$ the section $s$ lifts canonically to a section of $\pi_1(X'/k)$. 
The presence of a section implies that the index is a power of $p$, see \cite{stix:periodindex} Thm 15, so that $\Br(X'/k)$ is a cyclic $p$-group. In the limit over all finite subcovers $X'/X$ of $\tldxs$ of (\ref{eq:LerayGm}) we therefore find the exact sequence
\[
0 \to \Br(k) \otimes \bZ_\ell \to \Br(\tldxs) \otimes \bZ_\ell \to \varinjlim_{\tldxs/X'/X} \rH^1(k,\Pic_{X'}) \otimes \bZ_\ell.
\]
The degree sequence on every $X'$ gives an exact sequence
\[
0 \to \bZ/\pe(X')\bZ \to \rH^1(k,\Pic^0_{X'}) \to \rH^1(k,\Pic_{X'}) \to \rH^1(k,\bZ) = 0.
\]
The period, i.e., the order of $[\Pic^1] \in \rH^1(k,\Pic^0)$, divides the index and thus is a power of $p$ for any subcover $X' \to X$ of $\tldxs$. We therefore find in the limit
\[
\varinjlim_{\tldxs/X'/X} \rH^1(k,\Pic_{X'}) \otimes \bZ_\ell \cong \varinjlim_{\tldxs/X'/X} \rH^1(k,\Pic^0_{X'}) \otimes \bZ_\ell \cong  \rH^1(k, \varinjlim_{\tldxs/X'/X} \Pic^0_{X'}) \otimes \bZ_\ell \cong 0
\]
because $\varinjlim_{X_s/X'/X} \Pic_{X'}^0$ is uniquely divisible. Indeed, for every finite subcover $X'/X$ of $\tldxs$ and every $n \in \bN$ we have a further subcover $X'_n/X'$ of $\tldxs$, such that $\Pic_{X'}^0 \to \Pic_{X'_n}^0$ factors over the multiplication by $n$ map of $\Pic_{X'}^0$. Namely, if $X'$ corresponds to $H\cdot s(\Gal_k) \subset \pi_1(X)$ with $H \subset \pi_1(X \times_k k^\alg)$ open, then $X'_n$ corresponds to $\ov{[H,H]H^n} \cdot s(\Gal_k)$.
\end{pro}


\subsection{Detecting a valuation from a section} \label{sec:detectvaluation}
\begin{thm} \label{thm:createaplace}
Let $s : \Gal_k \to \pi_1(X)$ be a section. There is a valuation $w \in \Val_v(K)$ with prolongation  $\tilde{w}$ to $\tld K$ such that 
\begin{enumer}
\item the image  $s(\Gal_{k,\ell})$ of an $\ell$-Sylow subgroup $\Gal_{k,\ell}$ of $\Gal_k$  is contained in $\rD_{\tilde{w}|w}$, and
\item the image  $s(\rI_{k,\ell})$ of the $\ell$-Sylow subgroup $\rI_{k,\ell} = \Gal_{k,\ell} \cap \rI_k$ is contained in $\rI_{\tilde{w}|w}$.
\end{enumer}
\end{thm}
\begin{pro}
Let $M = k(\tldxs)$ be the function field of the decomposition pro-cover $\tldxs$ of the section $s$. 
By Theorem~\ref{thm:popdiss} and Proposition~\ref{prop:brauerxs} there is a valuation $w \in \Val_v(M)$, such that the Brauer class of invariant $1/\ell$ in $\Br(k)$ is nontrivial in the Brauer group of the henselisation $M^{\rh}_{w}$. We claim, that  $s(\Gal_{k,\ell})$ is contained in $\rD_{\tilde{w}|w}$ for an appropriate prolongation $\tilde{w} \in \Val_v(\tilde{K})$ of $w$.

With  $\Lambda=\tilde{K} \cap M^{\rh}_{w}$ the claim is equivalent to the degree of $\Lambda/M$ being prime to $\ell$ in the sense of supernatural numbers. By construction of $M$ there is an algebraic extension $\lambda/k$ such that $\Lambda=\lambda M$ and the degree of $\Lambda/M$ equals the degree of $\lambda/k$. The defining property of $w$ implies that $\Br(k) \otimes \bZ_\ell \inj \Br(\lambda) \otimes \bZ_\ell$ is injective, which forces the degree of $\lambda/k$ to be prime to $\ell$. This proves the claim and we have found a valuation $w$ such that (i) holds.

In order to enforce property (ii), let us first assume that the valuation $w$ constructed above is of \valht $2$. In the commutative diagram
\[
\xymatrix@M+1ex@R-2ex{ 1 \ar[r] & \rI_{\tilde{w}|w} \ar[d] \ar[r] & \rD_{\tilde{w}|w} \ar[d] \ar[r] & \Gal_{\kappa(w)} \ar@{^(->}[d] \ar[r] & 1 \\
 1 \ar[r] & \rI_k \ar[r] & \Gal_k \ar[r] & \Gal_\bF \ar[r] & 1 }
\]
the rightmost vertical map injective. Hence the $\ell$-Sylow subgroup $s(\rI_{k,\ell}) = s(\Gal_{k,\ell}) \cap \rI_k$ of $s(\rI_k)$ is automatically contained in $\rI_{\tilde{w}|w}$. 

It remains to discuss the case where a priori $w=\alpha$ is a valuation of $M$ of type 1v and $A$,  the Brauer class of invariant $1/\ell$ in $\Br(k)$, vanishes in $\Br(M_{w_{y\alpha}}^{\rh})$ for all valuations $w_{y\alpha} = v_y \circ \alpha$ of \valht $2$. 
The exact sequence (\ref{eq:valsequence}) yields in the limit an exact sequence
\begin{equation} \label{eq:valsequence2}
0 \to \Br(\kappa(\alpha)) \to \Br'(M_{\alpha}^{\rh}) \to \rH^1\big(\kappa(\alpha),\alpha(M) \otimes \bQ/\bZ\big) .
\end{equation}
The local--global principle for the Brauer group of $\kappa(\alpha)$ yields an injection
\[
\Br(\kappa(\alpha))  \inj \prod_{y}   \Hom\big(\Gal_{\kappa(w_{y\alpha})},  \bQ/\bZ\big) \inj  
 \prod_{y}  \Br\big(M_{w_{y\alpha}}^{\rh}\big)
\]
where $y$ ranges over the closed points of $Y_\alpha$ and so the composite valuations $w_{y\alpha} = v_y \circ \alpha$ are the refinements of type 2v of $\alpha$. Hence the restriction of $A$ in   $\Br'(M_{\alpha}^{\rh})$ is a ramified class, i.e., it induces a nontrivial character $\chi_A \in  \rH^1\big(\kappa(\alpha), \alpha(M) \otimes \bQ/\bZ\big)$. In fact $\chi_A$ is the character 
\[
\chi_A : \Gal_{\kappa(\alpha)} \to \Gal_\bF \xrightarrow{\Frob \mapsto v(\pi) \otimes 1/\ell} v(k) \otimes  \bQ/\bZ \to  \alpha(M) \otimes \bQ/\bZ,
\] 
where $\pi$ is a uniformizer of $\fok$. For $\chi_A$ to be nontrivial means that
the image of 
\[
\rD_{\tilde{\alpha}|\alpha} = \Gal(\tilde{K}/M_{\alpha}^{\rh}) \subset \Gal(\tilde{K}/M) = s(\Gal_k) \to \Gal_k \to \Gal_\bF
\]
contains the $\ell$-Sylow subgroup of $\Gal_\bF$, and the ramification index $e(\alpha/v)$ is prime to $\ell$.

Let $L/M$ be a subextension of $\tilde{K}/M$ with an unramified  prolongation $\alpha_L$ of $\alpha$ and such that the residue field of $\alpha_L$ has Galois group 
\[
\Gal(\kappa(\tilde{\alpha})/\kappa(\alpha_L)) \subset \Gal(\kappa(\tilde{\alpha})/\kappa(\alpha))
\]
that projects isomorphically to the  $\ell$-Sylow subgroup $\Gal_{\bF,\ell}$ of $\Gal_\bF$. Consequently, the restriction of $A$ to $\Br(L_{\alpha_L}^{\rh})$ still does not vanish because the restriction of the character $\chi_A$  remains nontrivial. By construction we have a short exact sequence
\begin{equation} \label{eq:decalphaL}
1 \to \rI_{\tilde{\alpha}|\alpha} \to \rD_{\tilde{\alpha}|\alpha_L} \to \Gal_{\bF,\ell} \to 1.
\end{equation}
The argument in the first part of the proof shows that (i) holds for $L$ and $\alpha_L$, thus showing that $s(\Gal_{k,\ell}) \subset \rD_{\tilde{\alpha}|\alpha_L}$. By a diagram chase with (\ref{eq:decalphaL}) we deduce that 
\[
s(\rI_{k,\ell}) \subset  \rI_{\tilde{\alpha}|\alpha} 
\]
and so $\alpha$ indeed also satisfies (ii).
\end{pro}


\section{Proof of the main result} \label{sec:fixedpt} 
In this section we formulate and prove a slightly stronger form of part (1) of the Main Result from the introduction.


\subsection{Reformulation as a fixed point problem}
The $T^G$ be the set of fixed points for a continuous $G$-action on a space $T$.


\subsubsection{The action}
The fundamental group $\pi_1(X,\bar{\eta})$ is anti-isomorphic to the group of covering transformations of the universal pro-\'etale cover $\tilde{X}/X$. Thus $\pi_1(X,\bar{\eta})$ acts on $\tilde{X}$ from the right and on its function field $\tilde{K}$ from the left. 
The action on $\tilde{X}$ is continuous in the sense that the induced action on a finite intermediate Galois cover $X'/X$ factors through a finite quotient of $\pi_1(X,\bar{\eta})$. Cofinally in the set of all intermediate covers $X'/X$ and their models $\cX'$ we find Galois equivariant models to which the $\pi_1(X,\bar{\eta})$-action uniquely extends.
The set of valuations $\Val_v(\tilde{K})$ of $\tilde{K}$ extending $v$ inherits a continuous $\pi_1(X,\bar{\eta})$-action from the right as the pro-finite limit of  $\pi_1(X,\bar{\eta})$-spaces
\[
\varprojlim_{X' \subset \cX'} \cX'_{\bF,\cons},
\]
where $\cX'$ ranges over a cofinal system of Galois equivariant models.


\subsubsection{Fixed points and decomposition subgroups} The stabilizer of a valuation $\tilde{w} \in \Val_v(\tilde{K})$ is nothing but the decomposition subgroup $\rD_{\tilde{w}|w} \subset \pi_1(X,\bar{\eta})$. Hence for a subgroup $G \subset \pi_1(X,\bar{\eta})$ we have $G \subseteq \rD_{\tilde{w}|w}$ if and only if $\tilde{w}$ belongs to the fixed points  $\Val_v(\tilde{K})^G$ of the induced $G$-action.


\begin{thm} \label{thm:sectionsarelocal}
Let $X$ be a smooth, hyperbolic, geometrically connected curve over a finite extension $k$ of $\bQ_p$. Then for any section $s$ of $\pi_1(X/k)$ there exists a valuation $\tilde{w} \in \Val_v(\tilde{K})$, such that the image $s(\Gal_k)$ is contained in the decomposition subgroup $\rD_{\tilde{w}|w}$.
\end{thm}

\begin{pro}
Let $\Sigma$ be the image of $s(\Gal_k)$. By the above we have to show that the set of fixed points 
\[
\Big(\Val_v(\tilde{K})\Big)^\Sigma = \Big(\varprojlim_{X' \subset \cX'} \cX'_{\bF,\cons}\Big)^{\Sigma}
= \varprojlim_{X' \subset \cX'} \big(\cX'_{\bF,\cons}\big)^{\Sigma}
\] 
is non-empty, where $X' \subset \cX'$ ranges over Galois equivariant models of the smooth projective compactifications of finite Galois \'etale covers $X'/X$ in $\tilde{X}$. But for each Galois equivariant model the set of fixed points $\big(\cX'_{\bF,\cons}\big)^{\Sigma}$ is a closed subset of the pro-finite set $\cX'_{\bF,\cons}$ and thus compact. The projective limit of compact sets is non-empty if and only if each member of the limit is non-empty, which reduces the proof to the following Theorem~\ref{thm:havefixedpoints}. 
\end{pro}


\begin{thm} \label{thm:havefixedpoints}
Let $\Sigma \subset \pi_1(X,\bar{\eta})$ be the image of a section $s : \Gal_k \to \pi_1(X,\bar{\eta})$, and let $X'/X$ be a finite Galois \'etale cover with a Galois equivariant model $\cX'$. The induced action of $\Sigma$ on $\cX'_{\bF,\cons}$ has a nonempty set of fixed points $(\cX'_{\bF,\cons})^\Sigma$.
\end{thm}
 
The proof of Theorem~\ref{thm:havefixedpoints} will be given in Section~\ref{sec:prooffixedpoint} below after recalling the following preliminary lemma.


\subsubsection{Sturdy reduction}
We recall the following result from \cite{Mz:globalprofinite} Lemma 2.9.

\begin{lem}[Mochizuki]  \label{lem:sturdy}
Every model $\cX$ admits a finite log \'etale cover $\cX' \to \cX$ such that every strict transform of a component of the stable model in $\cX'$ has genus at least $2$.
\hfill $\square$
\end{lem}

\begin{rmk}
(1)
A cover $\cX'$ with degeneration of its stable model as in Lemma~\ref{lem:sturdy} is called \textbf{sturdy} in \cite{Mz:globalprofinite}.

(2)
The cover $\cX'$ im Lemma~\ref{lem:sturdy} is usually not regular, but can be turned into a model by a minimal desingularisation of rational $A_n$-singularities in the nodes.
\end{rmk}


\subsection{The existence of fixed points: proof of Theorem~\ref{thm:havefixedpoints}} \label{sec:prooffixedpoint} 

For fine enough finite \'etale covers $X'/X$ the smooth compactification of $X'$ will itself be hyperbolic. It thus suffices to consider the case of smooth projective curves $X/k$ of genus at least $2$.

Let $\Sigma \subset \pi_1(X,\bar{\eta})$ be the image of a section $s : \Gal_k \to \pi_1(X,\bar{\eta})$, and let $\Theta \subset \Sigma$ be the image $s(\rI_k)$ of the inertia subgroup under the section $s$. 
Let $w \in \Val_v(K)$ be a valuation as in Theorem~\ref{thm:createaplace} with a prolongation $\tilde{w}  \in \Val_v(\tilde{K})$ such that an $\ell$-Sylow $\Sigma_\ell$ of $\Sigma$ is contained in the decomposition group $\rD_{\tilde{w}|w}$ and  the $\ell$-Sylow subgroup $\Theta_\ell = \Sigma_\ell \cap \Theta$ of $\Theta$ is contained in $\rI_{\tilde{w}|w}$. 

 Let $X'/X$ be a finite Galois \'etale cover with a Galois equivariant model $\cX'$.  In order to prove Theorem~\ref{thm:havefixedpoints} we may pass to a characteristic cover  $X''/X'$ and a finer equivariant model $\cX'' \to \cX'$. Hence we may assume that $X'$ has field of constants a finite extension $k'$ and  admits a stable model over the valuation ring $\fo'$ of $k'$. 
  Moreover, by Lemma~\ref{lem:disentangle} and Lemma~\ref{lem:sturdy} we may assume that 
 \begin{enumerate}
 \item[(i)] the stable model $\cX'_{\rm stable}$ of $X'$ is sturdy, i.e., that any stable component is of genus at least $2$,  and 
 \item[(ii)] any two components of the stable model intersect at most once in the stable model.
 \end{enumerate}
Let now $w'$ be the restriction of $\tilde{w}$ to $K'$, i.e., the extension of $w$ to $K'$ determined by $\tilde{w}$.   The intersection $\Theta_\ell \cap \rI_{\tilde{w}|w'}$ is of finite index in $\Theta_\ell$ and thus isomorphic to $\bZ_\ell(1)$. 

Let $\cN_{X' | \cX'}$ be the kernel of the specialisation map $\spez : \pi_1(X') \to \pi_1(\cX')$, which contains $\rI_{\tilde{w}|w'}$.  
The projection $\pi_1(X') \to \Gal_k$ induces a map $\cN_{X'|\cX'} \to \rI_k^{\tame} = \hat{\bZ}'(1)$, which maps $\Theta_\ell \cap  \rI_{\tilde{w}|w'}$ to a pro-$\ell$-group of finite index in the $\ell$-Sylow subgroup $\bZ_\ell(1)$ of $\rI_k^\tame$, hence to an infinite pro-$\ell$-group.
Consequently, the image of $\Theta_\ell \cap \rI_{\tilde{w}|w'}$ in  $\cN_{X' | \cX'}^\ab \otimes \bZ_\ell$ is nontrivial and also isomorphic to $\bZ_\ell(1)$.
 
For an element $\sigma \in \Sigma$ the image of $\Theta_\ell \cap \rI_{\tilde{w}|w'}$ in $\cN_{X' | \cX'}^\ab \otimes \bZ_\ell$ meets the image of its $\sigma$-conjugate $\sigma\Theta_\ell\sigma^{-1} \cap \rI_{\sigma(\tilde{w})|\sigma(w')}$ nontrivially because both are contained in the image of $\Theta \cap \cN_{X' | \cX'}$ and we have the following well known and useful lemma. 

\begin{lem}
Let $H \subset \rI_k$ be a closed subgroup of the inertia group $\rI_k$ of $k$. Then the maximal pro-$\ell$ quotient $H^\ell$ of $H$ for an $\ell \not=p$ is a quotient of $\bZ_\ell(1)$.
\end{lem}
\begin{pro}
The wild inertia $\rP_k \lhd \rI_k$ is the unique normal $p$-Sylow subgroup of $\rI_k$. Thus $H^\ell$ is a quotient of $H/(H \cap \rP_k)$ which is a subgroup of $\rI_k/\rP_k \cong \prod_{\ell \not= p} \bZ_\ell(1)$ .
\end{pro}
 
We deduce that the images of $\rI_{\tilde{w}|w'}^{\ab} \otimes \bZ_\ell$ and $\rI_{\sigma(\tilde{w})|\sigma(w')}^{\ab} \otimes \bZ_\ell$ in $\cN_{X' | \cX'}^\ab \otimes \bZ_\ell$ intersect nontrivially.
Due to Proposition~\ref{prop:compare} we may compute in the subquotient $\cN_{X' | \cX'}^{\log,\ab} \otimes \bZ_\ell$ of $\pi_1^{\log}(\cX')$.  If $w$ has \valht $2$, Corollary~\ref{cor:inertiawrtstablemodel} implies that the intersection 
\[
\Big(\bigoplus_{\alpha} \rI_{\tilde{\alpha}|\alpha}^\ab \otimes \bZ_\ell \Big) \cap \Big( \bigoplus_{\alpha} \rI_{\sigma(\tilde{\alpha})|\sigma(\alpha)}^\ab \otimes \bZ_\ell \Big)
\]
in $\cN_{X' | \cX'}^\ab \otimes \bZ_\ell$ is nontrivial, where $\alpha$ ranges over valuations of type 1v associated to irreducible components $Y_\alpha$ of the special fibre of $\cX'_{\rm stable}$ that contain the image 
$y$ of the center $x_{w'} \in \cX'_{\bF}$ of $w'$,  and the $\tilde{\alpha}$ are 
prolongations to $\tilde{K}$ as in Corollary~\ref{cor:inertiawrtstablemodel}. If $w$ is of type 1v, the same conclusion holds with just $\alpha=w$.
Now Proposition~\ref{prop:detect1v} applies and shows that the intersection
\[
\{\alpha \ ; \ y \in Y_\alpha\} \cap \{\sigma(\alpha) \ ; \ y \in Y_\alpha\}
\]
is nonempty for every $\sigma \in \Sigma$. If the set $\{\alpha \ ; \ y \in Y_\alpha\}$ has cardinality $1$ then this $\alpha$ is a fixed point. Otherwise the cardinality is $2$ and the combinatorics of the $\Sigma$-action on $\Sigma \cdot \{\alpha \ ; \ y \in Y_\alpha\}$ conforms to the following combinatorial lemma.

\begin{lem} \label{lem:combinatorial}
Let $G$ be a finite group acting on a set $M$. Let $x,y \in M$ be elements, such that $M = G\cdot x \cup G \cdot y$ and such that for every $g \in G$ the set $\{x,y\}$ intersects $\{gx,gy\}$ nontrivially, then we have one of the following three cases.
\begin{enumera}
\item $M^G \not= \emptyset$, more precisely $x$ or $y$ is fixed under $G$.
\item $M = \{x,y\}$ has two elements and $G$ acts transitively.
\item $M = \{x,y,z\}$ has three elements and $G$ acts transitively.
\end{enumera}
\end{lem}
\begin{pro}
For $z \in M$ let $G_z$ be the stabilizer of $z$ under the action by $G$. If $G$ acts with two orbits, then $G = G_x \cup G_y$ and thus not both stabilizers are of index in $G$ bigger than $1$, hence we have case (1). The same conclusion holds if $x$ equals $y$.

If $G$ acts transitively on $M$ and $x \not= y$, then there is a $g \in G$ with $gx=y$ and we have 
\[
G = G_x \cup gG_x g^{-1} \cup gG_x \cup G_xg^{-1}.
\]
Because $G_x \cap gG_x g^{-1}$ contains $1 \in G$, we can estimate $\# G +1 \leq 4 \cdot \# G_x$ and thus the index $(G:G_x)$ is at most $3$. This proves the lemma.
\end{pro}

In the situation of the proof of Theorem~\ref{thm:havefixedpoints}, when we let $\Sigma$ act through a finite quotient on $\Sigma \cdot \{\alpha \ ; \ y \in Y_\alpha\}$, then we obviously have a fixed point in case (1), namely the generic point of the component $Y_\alpha$ of the special fibre of $\cX'_{\rm stable}$. Consequently, the generic point of the strict transform of $Y_\alpha$ in $\cX'$ is fixed by $\Sigma$.

Next we lead case (3) to a contradiction. In case (3) the set $\Sigma \cdot \{\alpha \ ; \ y \in Y_\alpha\}$ consists of three distinct valuations $\alpha_1,\alpha_2,\alpha_3$ such that by Proposition~\ref{prop:detect1v} we have
\[
\bZ_\ell(1) \oplus \bZ_\ell(1) \oplus \bZ_\ell(1) = \bigoplus_{i=1}^3 \rI_{\tilde{\alpha}_i|\alpha_i}^\ab \otimes \bZ_\ell \subset \cN_{X'|\cX'}^\ab \otimes \bZ_\ell.
\]
Moreover, the nontrivial  images of the conjugates of $\Theta_\ell \cap \rI_{\tilde{w}|w'}$ in $\cN_{X'|\cX'}^\ab \otimes \bZ_\ell$ all agree. But on the other hand, if $\{\alpha \ ; \ y \in Y_\alpha\} = \{\alpha_1,\alpha_2\}$ then 
the image of $\sigma\Theta_\ell\sigma^{-1} \cap \rI_{\sigma(\tilde{w})|\sigma(w')}$ is contained in 
$\rI^\ab_{\sigma(\tilde{\alpha}_1)|\alpha_1} \oplus \rI^\ab_{\sigma(\tilde{\alpha}_2)|\alpha_2}$, so lies in a coordinate plane. Because $\Sigma$ acts transitively, we get a contradiction.  

In case (2) of the lemma we find at least a fixed point in the stable model $\cX'_{\rm stable}$ that is the unique node $y$, due to condition (ii), in which the two components $Y_\alpha$ with $y \in Y_\alpha$ meet. The proof of Theorem~\ref{thm:havefixedpoints} will thus be completed by the following lemma.

\begin{lem}
Let $X/k$ be a smooth projective curve of genus at least $2$ that admits a stable model $\cX_{\rm stable}$. Let $\cX$ be a model of $X$ that allows a finite group action by $G$.  The natural map $f: \cX \to \cX_{\rm stable}$ to the stable model is then $G$-equivariant and the map on fixed points
\[
\cX^G \to \cX_{\rm stable}^G
\]
is surjective.
\end{lem}
\begin{pro}
The uniqueness of the stable model induces an action of $G$ and forces the map $f$ to be $G$-equivariant. 

Let $y \in \cX_{\rm stable}$ be a fixed point under $G$. Then $f^{-1}(y)$ is geometrically connected and consists either of just one point, which then necessarily is fixed by $G$, or is a tree of projective lines. Then the dual graph $\Gamma_y$ of $f^{-1}(y)$ is a tree which inherits an action by $G$.  By Lemma~\ref{lem:tree} below, we have a fixed vertex or a fixed edge in $\Gamma_y$. That translates into a fixed component or a fixed node in $f^{-1}(y)$, so anyway the set of fixed points in $\cX$ above $y$ is nonempty. 
\end{pro} 

\begin{lem} \label{lem:tree}
Let $G$ be a group acting on a finite nonempty graph $\Gamma$. If $\Gamma$ is a tree, then the $G$-action on $\Gamma$ has a fixed point, which can be a vertex or an edge.
\end{lem}
\begin{pro} This follows at once from \cite{serre:trees} Prop 10 and its corollary, which unfortunately is only stated for trees of odd diameter, when the guaranteed fixed point is a vertex. We recall the argument in order to cover the case of even diameter.

For two vertices $x,y \in \Gamma$ the \textbf{distance} $d(x,y)$ is defined as the minimum over the number of edges in a connected subgraph of $\Gamma$ that contains $x$ and $y$, see \cite{serre:trees} \S I.2.2. We set 
\[
d_x = \max \{d(x,y) \ ; \ \text{ all vertices } y \in \Gamma\}
\]
for any vertex $x \in \Gamma$, and call $d = \max_x \{d_x\}$ the \textbf{diameter} of $\Gamma$, see \cite{serre:trees} \S I.2.2. 
The function $d_x$ is convex along \textbf{geodesic} paths in $\Gamma$, as can be seen from an easy case by case proof of $d_x + d_z \geq 2d_y$ for adjacent vertices
\[
\xymatrix@R=-3pt{{x} & {y} & {z}  \\
  {\bullet}   \ar@{-}[r] & {\bullet} \ar@{-}[r] & {\bullet}} 
\]
Let $\Gamma' \subset \Gamma$ be the minimal connected subgraph of $\Gamma$ that contains all vertices $x \in \Gamma$ such that $d_x < d$. By the convexity of the function $d_x$ along geodesics, we find that $\Gamma'$ does not contain vertices $x$ with $d_x = d$. The tree $\Gamma'$ is thus a $G$-equivariant subtree and has smaller diameter than $\Gamma$. In fact, the diameter of $\Gamma'$ is $d-2$. By induction on the diameter it suffices to treat cases, where $\Gamma'$ is empty. This leaves only the case of diameter $1$ and $2$ which are trivial.
\end{pro}


\section{Arithmetic properties of the valuations given by sections} \label{sec:exotic}

In this section we would like to discuss arithmetic properties of the valuations defined by sections as given by Theorem~\ref{thm:sectionsarelocal}. 


\subsection{Sections localized at type 2h valuations and the $p$-adic section conjecture}
Let  $s:\Gal_k \to \pi_1(X)$ be a section with $s(\Gal_k)$ contained in the decomposition group $\rD_{\tilde{w}|w}$ of a valuation $\tilde{w} \in \Val_v(\tilde{K})$ of type 2h. The valuation $w$ is a refinement of a valuation $w_a$ of type 1h corresponding to a closed point $a \in X$ of the generic fibre. It follows that 
\[
s(\Gal_k) \subseteq \rD_{\tilde{w}|w} = \rD_{\tilde{w}_a|w_a} = s_a(\Gal_{\kappa(a)}) \]
which after projection to $\Gal_k$ implies $\Gal_k \subseteq \Gal_{\kappa(a)} \subseteq \Gal_k$, and $s=s_a$ is the section associated to the $k$-rational point $a \in X(k)$ as predicted by the $p$-adic section conjecture.

The $p$-adic section conjecture thus reduces to the task of eliminating valuations $w \in \Val_v(K)$ of type other than 2h in Theorem~\ref{thm:sectionsarelocal}.


\subsection{The residue field}
Let $s:\Gal(k) \to \pi_1(X)$ be a section and let $\tilde{w} \in \Val_v(\tilde{K})$ such that with $s(\Gal_k) \subset \rD_{\tilde{w}|w}$. The induced map
\[
\Gal_k \to \rD_{\tilde{w}|w}/\rI_{\tilde{w}|w} \to \Gal_\bF
\]
is surjective. Hence the residue field $\bF$ of $k$ is relatively algebraically closed in $\kappa(w)$. Therefore, if  $w=\alpha$ is of type 1v, we find that $\kappa(\alpha)$ is a regular function field over $\bF$. And if $w$ has \valht $2$, then $\kappa(w)$ equals $\bF$.  We conclude that the valuation $w$ given by Theorem ~\ref{thm:sectionsarelocal} cannot be of  type 2u$_{\rm smooth}$.


\subsection{Sections localized at valuations of \valht \texorpdfstring{$2$}{2}}
Let $s:\Gal(k) \to \pi_1(X)$ be a section and let $\tilde{w} \in \Val_v(\tilde{K})$ be a valuation of \valht $2$ with $s(\Gal_k) \subset \rD_{\tilde{w}|w}$.

\subsubsection{The ramification} Let $\cX/\fok$ be a model of $X$ with reduced geometric special fibre $\ov{Y}$. The \textbf{ramification} of a section $s: \Gal_k \to \pi_1(X)$ with respect to the model $\cX$ is defined as the homomorphism
\[
{\rm ram}(s) = \spez \circ s |_{\rI_k} : \rI_k \to \pi_1(\ov{Y})
\]
of the composite of the restriction to the inertia subgroup $\rI_k \subset \Gal_k$ with the specialisation map $\spez: \pi_1(X) \surj \pi_1(Y)$. A section $s$ is called \textbf{unramified} with rerspect to $\cX$  if ${\rm ram}(s)$ is the trivial homomorphism. A section associated to a $k$-rational point that extends to an $\fok$-rational point of the model, in particaular any such for proper $X/k$, is necessarily unramified.

Let $X/k$ be a proper, smooth hyperbolic curve. The diagram 
\begin{eqnarray} \label{eq:ram}
\xymatrix@M+1ex@R-2ex{ 1 \ar[r] & \rI_{\tilde{w}|w} \ar[d] \ar[r] & \rD_{\tilde{w}|w} \ar[d] \ar[r] & \Gal_{\kappa(w)} \ar@{^(->}[d] \ar[r] & 1 \\
 1 \ar[r] & \rI_k \ar[r] \ar@{.>}@/_3ex/[u]_s  & \Gal_k \ar@/_3ex/[u]_s  \ar[r] & \Gal_\bF \ar[r] & 1 }
\end{eqnarray}
shows that for any given model $\cX$ with reduced special fibre $Y$ the ramification ${\rm ram}(s)$
of the section $s$ vanishes. Moreover, the induced section of $\pi_1(Y/\bF)$ is the section associated to the $\bF$-rational point given by the center $x_w$ of the valuation $w$ on $Y \subset \cX$. Of course, this is predicted by the $p$-adic section conjecture, but in general this is not known for an arbitrary section.

\subsubsection{The non-vanishing locus of constant Brauer classes}  \label{sec:nonvanish} By the diagram (\ref{eq:ram}) above, the section induces a splitting of the projection $\rI_{\tilde{w}|w} \surj \rI_k$. By the computation of log inertia groups the section thus  yields a splitting of the map 
\[
\Hom\big(w(K^\ast),\hat{\bZ}'(1)\big)  = \rI_{\tilde{w}|w}^\tame \surj \rI_k^\tame = \Hom\big(v(k^\ast),\hat{\bZ}'(1)\big) 
\]
and this means that $v(k^\ast) \inj w(K^\ast)$ has no cotorsion prime-to-$p$. By Corollary~\ref{cor:criterion} it follows that for all $\ell \not= p$ the map
\[
\Br(k) \otimes \bZ_\ell \to \Br(K_w^{\rh}) \otimes \bZ_\ell
\]
is injective. 

\subsubsection{Independence of \texorpdfstring{$\ell\not= p$}{l different from p}}
Although a valuation for which the constant Brauer class of invariant $1/\ell$ does not vanish is the starting point in the proof of Theorem~\ref{thm:createaplace}, in the course of the proof of Theorem~\ref{thm:havefixedpoints}  no effort is taken to keep this property. It turns out that at least for valuations of \valht $2$ that satisfy the claim of Theorem~\ref{thm:sectionsarelocal} the non-vanishing of the constant Brauer class of invariant $1/\ell$ is automatic. Moreover, the potential dependence of the choice of the auxillary prime $\ell$ different from $p$ does not play a role in the end.


\section{Uniqueness properties of the valuations given by sections} \label{sec:unique}


\subsection{Bridges and the effect of resolution of non-singularities} \label{sec:bridgesnonres}

With regards to uniqueness of the valuation in Theorem~\ref{thm:sectionsarelocal}, we first discuss the combinatorial structure of the union of all invisible components. 


Because of Lemma~\ref{lem:visible} the invisible components are of genus $0$ over some field extension $\bF'$ of $\bF$ and meet the rest of the special fibre in at most a divisor of degree $2$ over $\bF'$.
A special kind of  invisible component is defined as follows.

\begin{defi} \label{defi:bridge}
A \textbf{bridge element} is an invisible irreducible component of the special fibre of a model $\cX$, which is contained in a \textbf{bridge}, i.e, a chain of components  $E_0 = Y_{\alpha}, E_1, \ldots,E_{e-1}$, $E_e = Y_{\alpha'}$ in the reduced special fibre $\cX_{\bF,\redu}$ where 
\begin{enumer}
\item $E_i$ meets $E_{i-1}$ and $E_{i+1}$ in a double point, 
\item $E_i$ is invisible for $i = 1, \ldots,e-1$,
\item $Y_{\alpha}$ and $Y_{\alpha'}$ are visible and not necessarily distinct components, the \textbf{bridge heads} of the bridge.
\end{enumer}
A valuation of type 1v is called a \textbf{bridge element} if the associated irreducible component on a fine enough model is a bridge element. Valuations $\alpha,\alpha'$ of type 1v which give rise to the bridge heads $Y_\alpha, Y_{\alpha'}$ are also called \textbf{bridge heads}.
\end{defi}

\begin{rmk}
(1) 
Due to Lemma~\ref{lem:visible}, a bridge element can only be dominated by bridge elements in refinements of models or in models of finite, generically \'etale covers. Hence a valuation of type 1v belongs to a bridge on every model on which its associated divisor appears.

(2) 
An unproven stronger form of \textit{resolution of non-singularities}, 
see \cite{tamagawa:resolution},  would imply that there are no invisible components and therefore also no bridges.
\end{rmk}


For the sake of reference we extract the following lemma from Tamagawa's work on resolution of non-singularities \cite{tamagawa:resolution}.

\begin{thm}[Tamagawa] \label{thm:resolution}
Let $y_1,y_2$ be distinct closed points on a visible component $Y_\alpha$ of the reduced special fibre of a model $\cX$ of $X/k$. Then there is a finite \'etale cover $X' \to X$ and a model $\cX'$ of $X'$ which allows an extension of the cover $f:\cX' \to \cX$, such that the following holds.
\begin{enumer}
\item We have distinct visible components $Y_{\alpha_1},Y_{\alpha'}, Y_{\alpha_2}$ in the reduced special fibre of $\cX'$.
\item $Y_{\alpha'}$ dominates $Y_\alpha$ under the map $f$.
\item  $f(Y_{\alpha_i}) = y_i$ for $i=1,2$.
\item  $Y_{\alpha'}$ and $Y_{\alpha_i}$ for $i = 1,2$ intersect above $y_i$ or are connected by a bridge, the bridge elements of which map to $y_i$ under $f$.
\end{enumer}
\end{thm}
\begin{pro}
That we can find a cover $X' \to X$ with a model $\cX'$  that satisfies (i)-(iii) follows directly from \cite{tamagawa:resolution} Thm 0.2 (v). Note that \cite{tamagawa:resolution} works with components of the stable model. By Lemma~\ref{lem:visible}, an auxiliary cover allows first to replace $Y_\alpha$ by a component of the stable model. 

Then, as our models are assumed to be regular, we have to resolve the singularities of the stable model,  that is rational $\rA_n$-singularities, which only contributes additional chains of $\bP^1$'s. By choosing $Y_{\alpha_i}$ visible and at minimal distance from $Y_{\alpha'}$ along such a chain yields the desired components.
\end{pro}


\subsection{Uniqueness of the valuation}
It is a natural question whether for a given section $s$ the valuation $\tilde{w} \in \Val_v(\tilde{K})$  given by  Theorem~\ref{thm:sectionsarelocal} such that $s(\Gal_k) \subseteq D_{\tilde{w}|w}$ is unique.

\begin{prop}  \label{prop:ht2}
Let $w_1,w_2$ be valuations of \valht $2$ with $s(\Gal_k) \subseteq \rD_{\tilde{w}_i|w_i}$ for $i=1,2$.
Let $\cX'$ be a Galois equivariant model of a finite \'etale Galois cover $X' \to X$. Then the centers $y'_i= x_{w'_i}$ of the $w'_i = \tilde{w}_i|_{X'}$ map to the same closed point in the stable model $\cX'_{\rm stable}$.
\end{prop}

\begin{pro}  
We use the notation of the proof of Theorem~\ref{thm:havefixedpoints}. By Lemma~\ref{lem:disentangle} and Lemma~\ref{lem:sturdy} we may assume that the stable model $\cX'_{\rm stable}$ of $X'$ is sturdy, i.e., that any stable component is of genus at least $2$,  and any two components of the stable model intersect at most once in the stable model. In order to simplify notation we assume that $X'=X$ with stable model $\cX_{\rm stable}$.

In $\cN_{X|\cX}^{\log,\ab} \otimes \bZ_\ell$ the log inertia groups $\rI_{w_i}^{\log} \otimes \bZ_\ell$ meet in the image $\Theta_\ell$ of an $\ell$-Sylow of $\rI_k$ under the section $s$. We argue as in the proof of Theorem~\ref{thm:havefixedpoints} using Proposition~\ref{prop:detect1v} that the irreducible components of the special fibre $Y$ of $\cX_{\rm stable}$ which contain $y_1=x_{w_1}$ cannot be disjoint from those which contain $y_2=x_{w_2}$. So there is at least one component $Y_\alpha$ corresponding to a valuation $\alpha$ of type 1v which contains both $x_{w_1}$ and $x_{w_2}$. 

We apply the preceding paragraph to finite \'etale covers $X'' \to X' \to X$ which are generic fibres of finite \'etale covers $\cX'' \to \cX'$. We deduce that the section of $\pi_1(Y/\bF)$ induced by $s$, namely $s_{y_1}=s_{y_2}$, factors as the corresponding section of $\pi_1(Y_\alpha/\bF)$. The injectivity of the natural map
\[
Y_\alpha(\bF) \to \left\{\text{conjugacy classes of sections of } \pi_1(Y_\alpha/\bF) \right\}
\]
shows thus that $y_1=y_2$ as claimed.
\end{pro}

\begin{prop} \label{prop:visible}
Let $s:\Gal_k \to \pi_1(X)$ be a section. Then there is at most one valuation $\alpha \in \Val_v(K)$  of \valht $1$ corresponding to a visible component $Y_\alpha$ of the special fibre of some model $\cX$ such that $s(\Gal_k) \subset \rD_{\tilde{\alpha}|\alpha}$.  Moreover, if such an $\alpha$ exists, then there is 
\begin{enumera}
\item either a refinement $w=v_y\circ \alpha$ of type 2v associated to a closed point $y\in Y_\alpha$ such that even $s(\Gal_k) \subset \rD_{\tilde{w}|w}\subset \rD_{\tilde{\alpha}|\alpha}$,
\item or the image $\Theta_\ell$ under $s$ of an $\ell$-Sylow of the inertia group $\rI_k \subset \Gal_k$ is contained in $\rI_{\tilde{\alpha}|\alpha}$ for all $\ell$.
\end{enumera}
\end{prop}
\begin{pro}
It follows essentially from Tamagawa's work on resolution of non-singularities \cite{tamagawa:resolution}, more concretely from the assertion of Theorem~\ref{thm:resolution}, that the Galois extension of the residue field  $\kappa(\alpha)$ at $\alpha$ corresponding to $\Gal_{\kappa(\alpha)} \surj \Gal_\alpha:=\rD_{\tilde{\alpha}|\alpha}/ \rI_{\tilde{\alpha}|\alpha}$ has no prime-to-$p$ extensions. Indeed, in the system of components $Y'_\alpha$ corresponding to $\tilde{\alpha}$ for finer and finer \'etale covers  $X'/X$ the set of nodes on $Y'_\alpha$ will contain any given set of closed points. But as 
\cite{Mz:globalprofinite} Prop 4.2 or  \cite{stix:diss} Prop 6.2.11, show that with the natural log structures  
\[
\pi_1^{\log}(Y'_\alpha) \inj \pi_1^{\log}(\cX',\bar{\eta}) 
\]
is injective, we see that any tamely ramified cover of $Y'_\alpha$ with ramification in the set of nodes will occur as residue extension of $\kappa(\alpha)$.

Let us now assume that $\Theta_\ell$ is not contained in $\rI_{\tilde{\alpha}|\alpha}$, so (2) fails. 
The proposition then claims, that for every Galois equivariant model $\cX'$ of a finite \'etale cover $X'/X$ there is a fixed point $y$ under the action of $\Sigma = s(\Gal_k)$ which is a closed point in the closure $Y'_\alpha$ of the center of $\alpha$. Moreover, there is a compatible system of such fixed points as the model varies. By Lemma~\ref{lem:disentangle} and Lemma~\ref{lem:sturdy} we may assume that the stable model $\cX'_{\rm stable}$ of $X'$ is sturdy, i.e., that any stable component is of genus at least $2$,  and any two components of the stable model intersect at most once in the stable model. Furthermore, it suffices to find such a fixed point $y$ on the stable model $\cX'_{\rm stable}$.

The image $\ov{\Theta}_\ell$ of $\Theta_\ell$ in $\Gal_w = \Gal(\kappa(\tld w)|\kappa(w))$ will be a cyclotomically normalized subgroup, see \cite{naka:1994} \S2.1, isomorphic to $\bZ_\ell(1)$ of $\ker(\Gal_w \surj \Gal_{\bF})$. The theory of the anabelian weight filtration as pioneered by Nakamura in \cite{naka:1990} \S3, \cite{naka:1994} \S2.1, see also \cite{stix:habil} \S26.6, still works in this context, because $\Gal_w$ is sufficiently big, and shows that $\ov{\Theta}_\ell$ is contained in an inertia subgroup of a unique node $y$ of a suitable  corresponding $Y'_\alpha$. 
By structure transport using conjugation by elements of $\Gal_k$ through the section $s$ we see that in fact the point $y$ is preserved under $\Sigma$.

By passing to finer and finer covers $X'$ and models $\cX'$ we deduce from the uniqueness of $y$ which is detected by the partial image $\ov{\Theta}_\ell$ of the section in $\Gal_w$ that the collection of closed points so obtained forms a compatible system in 
\[
\varprojlim_{X' \subset \cX'} \cX'_{{\rm stable},\bF}
\]
endowed with the constructible topology.
Because the valuation $\alpha$ corresponds to a visible component $Y'_\alpha \subset  \cX'_{\rm stable} $ that is  fixed by $\Sigma$, and because the system of closed points given by the $y$ lies on the $Y'_\alpha$, we may conclude that the corresponding closed points on the strict transforms of the $Y'_\alpha$  in any model are also preserved by $\Sigma$. Hence there is a refinement $w=v_y\circ \alpha$ of type 2v associated to the system of closed points $y\in Y'_\alpha$ such that $s(\Gal_k) \subset \rD_{\tilde{w}|w}\subset \rD_{\tilde{\alpha}|\alpha}$, as claimed by option (1).

It remains to prove the assertion on uniqueness. We can argue by Proposition~\ref{prop:detect1v} as in the proof of Proposition~\ref{prop:ht2}. Indeed, under option (1) or (2) the image $\Theta_\ell$ will detect the corresponding valuation of type 1v. The only problem that might occur is solved by moving apart the the two stable components $Y_\alpha, Y_\beta$ meeting in $y$ by an application of 
Theorem~\ref{thm:resolution}.
\end{pro}

\begin{prop} \label{prop:invisible}
Let $\cX'$ be a Galois equivariant model of a finite \'etale Galois cover $X' \to X$. 
Let $w_1,w_2$ be valuations  with $s(\Gal_k) \subseteq \rD_{\tilde{w}_i|w_i}$ for $i=1,2$, such that the  centers $x_{w'_i}$ of the $w'_i = \tilde{w}_i|_{X'}$ map to closed points $y'_i$  in the stable model $\cX'_{\rm stable}$. Then $y'_1$ coincides with $y'_2$.
\end{prop}
\begin{pro}
By assumption we have $\rD_{\tilde{w}_i|w'_i} \subset \rD_{y'_i}$ whose image in $\pi_1(Y')$ under the specialisation map coincides with the image of the sections associated to the closed points $y'_i$ in the reduced special fibre $Y'$ of the stable model $\cX'_{\rm stable}$. Replacing the $\rI_{w'_i}^{\log}$ by the $\rI_{y'_i}^{\log}$ now the proof of Proposition~\ref{prop:ht2} applies mutatis mutandis.
\end{pro}

\begin{thm} \label{thm:pseudouniqueness}
Let $\Sigma$ be the image of a section $s:\Gal(k) \to \pi_1(X)$. Let  $\cX'$ be a Galois-equivariant model of a finite \'etale cover $X' \to X$ with stable model $\cX'_{\rm stable}$. Then the image of the map
\[
{\rm center} : \Val_v(\tilde{K})^\Sigma \to \cX'_{{\rm stable},\bF}
\]
consists either
\begin{enumera}
\item of a unique closed $\bF$-rational point, or
\item of the generic point of a unique component together with a closed $\bF$-rational point on that component, or
\item of the generic point of a unique component. 
\end{enumera}
\end{thm}
\begin{pro}
Lemma~\ref{lem:visible}, Proposition~\ref{prop:visible} and Proposition~\ref{prop:invisible} show that the image contains at most one closed and at most one generic point, while Theorem~\ref{thm:sectionsarelocal} shows that the image is nonempty. It remains to argue that if the image consists of both a closed point $y$ and a generic point $\alpha$, then $y$ is in the closure of $\alpha$. But this follows clearly from Proposition~\ref{prop:detect1v} and the discussion of the group $\Theta_\ell$ in the proof of both Proposition~\ref{prop:visible} and Proposition~\ref{prop:invisible}.
\end{pro}

\begin{cor}
Let $X/k$ be a smooth hyperbolic curve which has a cofinal system of finite \'etale covers $X' \to X$ such that $X'$ has a model $\cX'$ whose components of the special fibre are visible. Let $s : \Gal_k \to \pi_1(X)$ be a section. Then one of the following holds. 
\begin{enumera}
\item  There is a unique $\tilde{w} \in \Val_v(\tilde{K})$ with  $s(\Gal_k) \subset \rD_{\tilde{w}|w}$.
\item  There exist a unique $\tilde{\alpha}$ of type 1v with a refinement $\tilde{w}=\tilde{v}_y \circ \tilde{\alpha}$  of type 2v, such that $s(\Gal_k) \subset \rD_{\tilde{w}|w} \subset \rD_{\tilde{\alpha}|\alpha}$.
\end{enumera}
\end{cor}
\begin{pro*}
This is an immediate corollary of Theorem~\ref{thm:pseudouniqueness} as the assumption of all components being visible leads to a bijection
\[
\ \hspace{4cm} \Val_v(\tilde{K}) \xrightarrow{\sim} \varprojlim_{X' \subset \cX'} \cX'_{ \bF} \xrightarrow{\sim} \varprojlim_{X' \subset \cX'} \cX'_{{\rm stable},\bF} \hspace{4cm} \qed
\]
\end{pro*}


\subsection{Final remark} 
It is conceivable that 
one may extend the range of uniqueness of Theorem~\ref{thm:pseudouniqueness} to also include the locus in bridges of suitable models. But as soon as there are invisible 
$\bP^1$'s for a curve $X/k$, it is also conceivable that those will contribute sections of $\pi_1(X/k)$ localized in the respective $p$-adic disc of the associated rigid analytic space but not localized in a $k$-rational point, thus ultimately failing the $p$-adic section conjecture. 

\begin{appendix}

\section{The zoo of valuations for  \texorpdfstring{$2$}{2}-dimensional semilocal fields} \label{sec:val}

Let $k$ be a complete discrete valued field with valuation ring
$\fok$ and perfect residue field $\kkp$, e.g., $k$ is a finite 
extension of $\bQ_p$. Let $v$ denote the canonical valuation 
on $k$.

\subsection{The Riemann--Zariski space of  $\fok$-valuations} 

Let $K$ be the function field of a smooth projective   
geometrically connected curve $X$ over $k$. In this appendix we 
discuss the space 
\[
\Val_\fok(K)=\Val_k(K)\cup\Val_v(K) =  \{w \ ; \  \text{ valuation on } K  \text{ with } w(\fo) \geq 0\} 
\]
of  equivalence classes of valuations $w$ on $K$ whose valuation
ring $R_w$ contains $\fok$, or equivalently, the restriction of $w$ to $k$ 
is either the trivial valuation or equals $v$.

\subsubsection{Models} \label{sec:model}
In this paper a \textbf{model} or more precisely a  \textbf{regular 
model with strict normal crossing} of $X$ over $\fok$ is a 
regular scheme $\cX$ which is flat and proper over 
$\Spec(\fok)$ together with a $k$-isomorphism of $X$ 
with the generic fibre $\cX_k$ such that the reduced 
special fibre $\cX_{\kkp,\redu}$ is a divisor with strict 
normal crossings on $\cX$.  In particular, unfortunately 
a stable model in general is not a model in the sense 
of this paper. By a result of Lichtenbaum, 
\cite{lichtenbaum:curvesdvr} Thm 2.8, models are 
automatically projective over $\Spec(\fok)$. We denote the 
underlying topological space of $\cX$ by $\cX^\topo$, 
whereas $\cX_\cons$ denotes $\cX^\topo$ when given 
the finer constructible topology.

\subsubsection{The center} \label{sec:center}
 For $w \in \Val_\fok(K)$ the 
valuative criterion of properness implies a canonical 
map $\Spec (R_w) \to \cX$ which maps the closed point 
of $\Spec(R_w)$ to the \textbf{center} $x_w \in \cX_\cons$ 
of the valuation $w$ on the model $\cX$. Maps between 
different models, which are the identity on $X$, respect 
the center of a valuation. The resulting map 
\begin{equation} \label{eq:center}
{\rm center} \, : \  \Val_\fok(K) \to \varprojlim \cX_\cons,
\end{equation}
where the projective limit ranges over all models of $K$, identifies 
$\varprojlim \cX_\cons$ with $\Val_\fok(K)$ which is a subspace
of the \textbf{Riemann--Zariski space} of $K/k$.

\subsubsection{The valuation ring} 
The centers of a valuation $w \in \Val_\fok(K)$  determine the valuation ring $R_w = \varinjlim \OO_{\cX,x_w}$ where the direct limit ranges over all models. The inverse map to (\ref{eq:center}) is described as follows. To a compatible system of points  $a_{\cX} \in \cX_\cons$ on all models we associate first the ring $R_a = \varinjlim \OO_{\cX\hhb{-2},\hhb1a_\cX}$. The ring $R_a$ is the valuation ring of a valuation $w$ of $K$ because for every $f \in K^\ast$ at least one of $f$ and $f^{-1}$ belongs to $R_a$, see \cite{bourbaki:comalg1-7} VI \S1.2 . Indeed, the indeterminacy of $f$, that is the set of points where neither $f$ nor $f^{-1}$ is defined, disappears on a fine enough model.

\subsubsection{The patch topology}
The \textbf{patch topology} on $\Val_\fok(K)$
is defined as the topology induced from the pro-finite product topology by the injective map
\[
{\rm sign} \, : \  \Val_\fok(K) \inj \prod_{f \in K^\ast} \{ - ,0,+\}
\]
that assigns to a valuation $w$ the collection of \textbf{signs} of the value $w(f)$ for each $f \in K^\ast$, where the sign of $f$ is $+$ if $w(f) >0$, it is $-$ if $w(f) < 0$ and the sign is $0$ if $w(f)=0$.
The condition on a collection of signs to belong to a valuation ring, namely that the subset in $K$ of $0$ and the nonnegative elements forms a ring which contains at least one of $f,f^{-1}$ for each $f \in K^\ast$, is a closed condition. Hence $\Val_\fok(K)$ is a pro-finite space, in particular it is  compact and Hausdorff.

The map ${\rm center}:   \Val_\fok(K) \to \varprojlim \cX_\cons$ defined in~(\ref{eq:center}) is a homeomorphism from $\Val_\fok(K)$ endowed with the patch topology to $\varprojlim \cX_\cons$ with respect to the $\varprojlim$-topology.
The subset
\[
\Val_v(K) = \{w \in \Val_\fok(K) \ ; \ w|_k=v\} \subset \Val_\fok(K)
\]
is a closed subset in the patch topology described by the condition that $w(\pi) > 0 $ for a uniformizer $\pi$ of $\fok$. The set $\Val_v(K)$ corresponds to the subset $\varprojlim \cX_{\kappa,\cons} \subset \varprojlim \cX_{\cons}$ where $\cX_{\kappa} \subset \cX$ is the special fibre.
 
\subsection{Types of valuations} 

We sketch the classification of the zoo of valuations and fix the terminology.

\subsubsection{\Valht} 
We define the  \textbf{\valht}of a  valuation $w \in \Val_\fok(K)$ as the well defined number 
$0,1$ or $2$ given by the height  $\height(x_w) = \dim(\OO_{\cX,x_w})$ in the sense of scheme theory of its center $x_w \in \cX$ 
 for all sufficiently fine models $\cX$ with respect to the system of all models. The unique valuation of \valht $0$ is the trivial valuation. 
 
\subsubsection{\Valht 1} \label{sec:ht1}
 Valuations of \valht $1$ are the discrete valuations associated to prime divisors on an arbitrary model fine enough such that the respective divisor appears. The corresponding prime divisor is either \textbf{vertical}, i.e., it maps to the closed point of $\Spec(\fok)$, or \textbf{horizontal}, i.e., it maps finitely to $\Spec(\fok)$. The first are called of \textbf{type 1v} whereas the latter valuations are called of \textbf{type 1h}. 
 
\begin{nota} \label{note:1v}
The usual notation for a valuation of $K$ of type 1v will  be $\alpha$. The corresponding prime divisor of a fine enough model will be denoted by $Y_\alpha$ and by abuse of notation has a generic point denoted by $\alpha$ again. The precise meaning of $\alpha$ will  always be clear from the context.
\end{nota}
 
\subsubsection{\Valht 2}
All the remaining valuations are of \valht $2$ and thus have all their centers at  
closed points of the special fibre. 
  
Let $w$ be a valuation of \valht $2$. For each valuation $\alpha$ of \valht $1$ we define the distance of $w$ to $\alpha$ on the model $\cX$ as the infimum of the number of irreducible components in a $1$-dimensional connected subscheme $Z \subset \cX$ which contains the center of $\alpha$ and of $w$. If $w$ keeps finite distance to any valuation of \valht $1$ as we vary over the system of all models, then there is a unique valuation $\alpha$ of \valht $1$ with a closed point $y$ on the associated divisor  such that $w$ is the composition of $\alpha$ with the valuation $v_y$ on the residue field of $\alpha$ associated to $y$. So $w = v_y \circ \alpha$ is called of \textbf{type 2v} (resp.\ 
\textbf{type 2h}) if $\alpha$ is \textbf{vertical} (resp.\ \textbf{horizontal}.)

The remaining valuations have centers which move away from any valuation of \valht $1$ and are called of \textbf{(coarse) type 2u (unbounded)}.

\subsection{Rigid analytic viewpoint}  \label{sec:rigview}
Valuations of \valht $2$ can be understood in terms of the associated rigid analytic space $X^\rig$. For a model $\cX$ we get a specialisation map 
\[ {\rm sp}_\cX: X^\rig \to \cX_{\kkp,\redu}\]
 from the rigid space to the set of closed points of the special fibre. 
 The preimage of a smooth closed point of $\cX_{\kkp,\redu}$ is 
 an open disc, the preimage of a node of $\cX_{\kkp,\redu}$ is an 
 annulus.

To a valuation $w$ of \valht $2$ we associate the system  $C_\cX = {\rm sp}_\cX^{-1}(x_w)$ of preimages of the centers indexed by the system of all models. The system of subsets $C_\cX$ is monotone decreasing with respect to inclusion when the model becomes finer. The valuation is uniquely determined by the system of the $C_\cX$ as
\[ R_w = \bigcup_{\cX} \OO_\cX(C_\cX) = \bigcup_{\cX} \big\{f \in K; \ f \text{ defined on } C_\cX, \|f\|_{C_\cX,\infty} \leq 1 \big\}, \]
where $\|f\|_{C_\cX,\infty}$ is the sup-norm of $f$ on $C_\cX$. The various types belong to distinctive geometric pictures of the system of the $C_\cX$ as follows.

\subsubsection{Type 2h}
For fine enough models, $C_\cX$ is an open disc with fixed center $x \in X^\rig$ and radius converging to $0$ with finer and finer models.

\subsubsection{Type 2v}
For fine enough models, $C_\cX$ is an  annulus, such that the correspponding annuli for finer and finer models share one common boundary. 

\begin{figure}[htbp]
  \centering
  \begin{minipage}[b]{6cm}
     \includegraphics*[height=5cm]{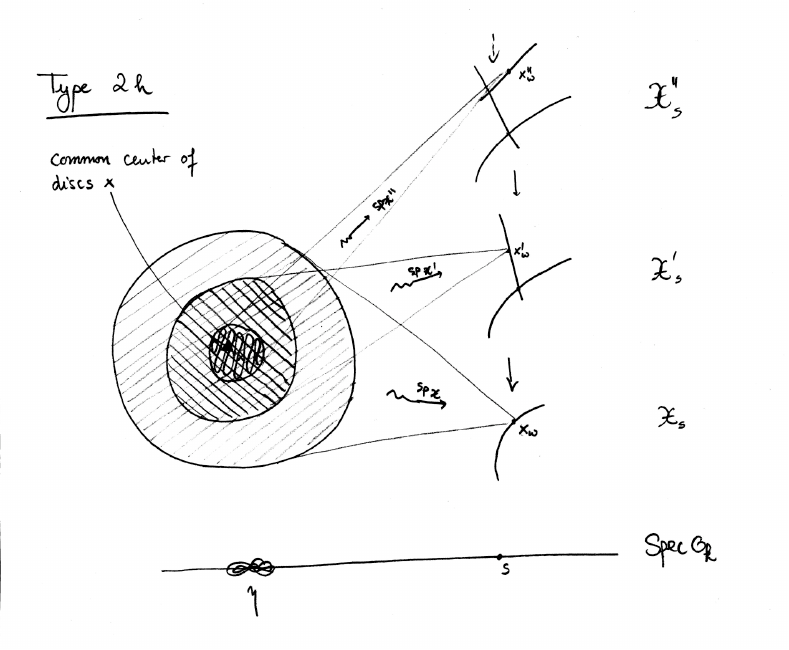}
     \caption{type 2h}
  \end{minipage}
  \hspace{1cm}
  \begin{minipage}[b]{6cm}
     \includegraphics*[height=5cm]{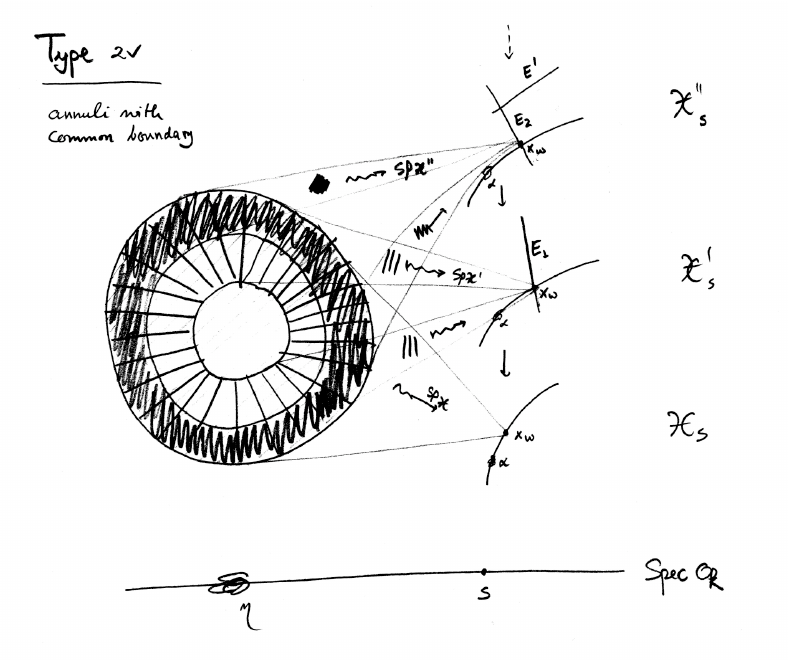}
     \caption{type 2v}
  \end{minipage}
\end{figure}

\subsection{\Valht 2 but unbounded distance}
The valuations of type 2u can be described and arranged into types in more detail as follows.

Every closed point $y$ in the reduced special fibre carries invariants $(e_{y,\alpha}),f_y$ equal to the tuple $(e_{y,\alpha})$ of the multiplicities of the components on which $y$ lies in the special fibre $\cX_\kappa$ and the residue field degree $f_y$ of $y$ over $\kkp$. Any closed  point $x$ in the generic fibre $X=\cX_K$ that specialises to $y$ has to have residue field $\kappa(x)$ with $e_{\kappa(x)/k} = \sum_\alpha m_\alpha e_{y,\alpha}$ with $m_i \in \bN_{\geq 1}$ and $f_y | f_{\kappa(x)/k}$. On the other hand, there is always an $x$ with the minimal possible values of $e,f$. 

\subsubsection{Type \texorpdfstring{2u$_{\rm smooth}$}{2u-smooth}}
For a  valuation $w$ of type 2u the value $\sum_\alpha e_{x_w,\alpha}$ remains bounded if and only if for fine enough models ultimately all centers $x_w$ belong to the smooth locus of the reduced special fibre. Such a valuation is called of \textbf{type 2u$_{\rm smooth}$} or \textbf{2u$_{\rm sm}$ (ultimately smooth)}.

\subsubsection{Type \texorpdfstring{2u$_{\rm node}$}{2u-node}}
We call a valuation $w$ of \textbf{type 2u$_{\rm node}$} or \textbf{2u$_{\rm n}$ (ultimately node)}  if for all fine enough models the center lies in a node of the reduced special fibre.

\subsubsection{Type \texorpdfstring{2u$_{\rm alt}$}{2u-alt}}
For a valuation of type 2u, if neither type 2u$_{\rm node}$ nor type 2u$_{\rm smooth}$ applies, then the center $x_w$ in the pro-system of models alternate between the smooth locus of the reduced special fibre and its nodes, and hence these are called of \textbf{type 2u$_{\rm alternating}$} or \textbf{2u$_{\rm alt}$(unbounded alternating)}.

\begin{figure}[htbp]
  \centering
  \begin{minipage}[b]{6cm}
     \includegraphics*[height=5cm]{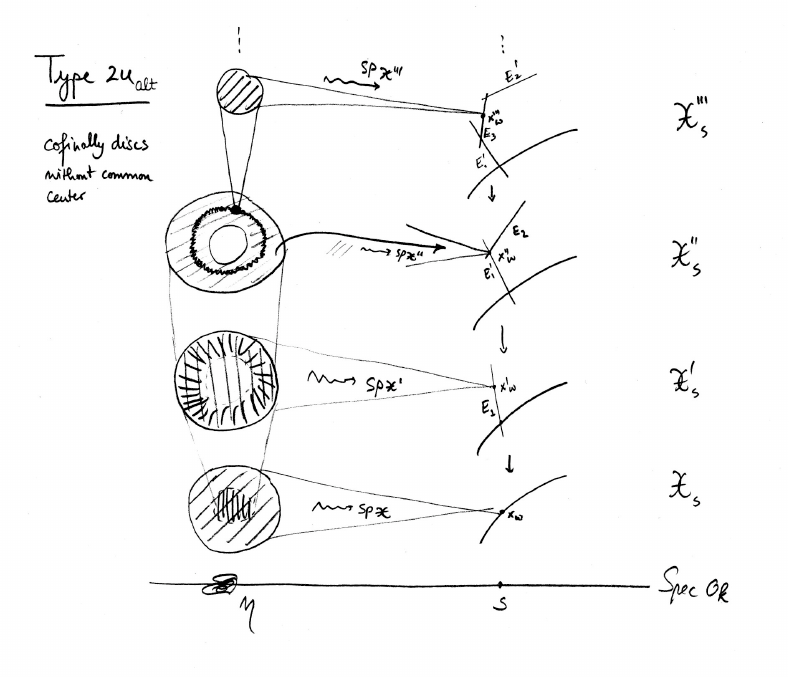}
     \caption{type 2u$_{\rm alternating}$}
  \end{minipage}
  \hspace{1cm}
  \begin{minipage}[b]{6cm}
     \includegraphics*[height=5cm]{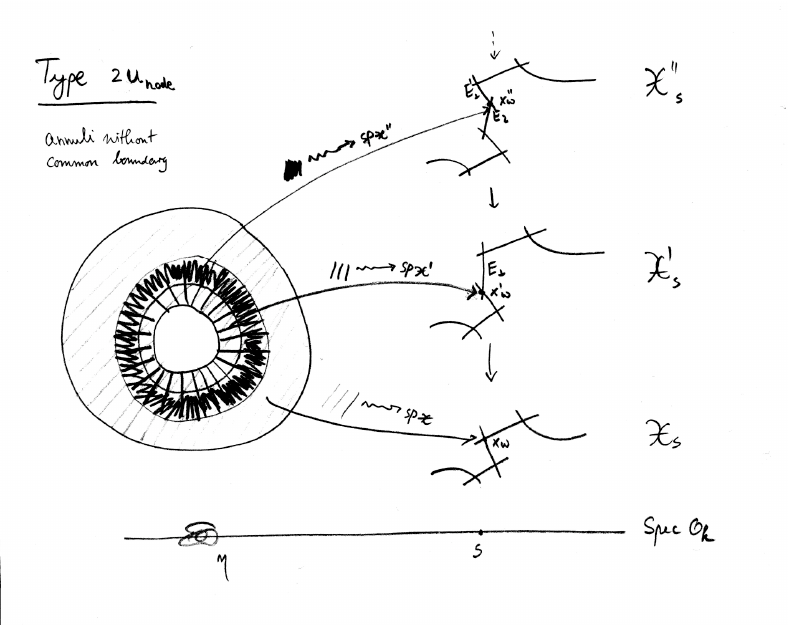}
     \caption{type 2u$_{\rm node}$}
  \end{minipage}
\end{figure}

\subsubsection{Rigid analytic description of type \texorpdfstring{2u$_{\rm smooth}$}{2usmooth}}
For a cofinal set of models, $C_\cX$ is an open disc without common center in $X^\rig$. The radius of the discs converges to $0$ with finer and finer models. There is a unique limit point in $X(\widehat{k^\alg}) \setminus X(k^\alg)$, where $\widehat{k^\alg}$ is the completion of $k^\alg$. 

\subsubsection{Rigid analytic description of type \texorpdfstring{2u$_{\rm node}$}{2unode}}
For fine enough models, $C_\cX$ is a $p$-adic  annulus, such that the corresponding annuli for finer and finer models share no common boundaries.

\subsection{Algebraic structure}
The information on the algebraic structure associated to a valuation $w$ according to its type is summarized in the following table. The \textbf{rational rank} of $w$ or better its value group $\Gamma_w$ is defined as $\dim_\bQ(\Gamma_w \otimes \bQ)$, see \cite{bourbaki:comalg1-7} VI \S10.2. And the \textbf{rank} of $w$,  \textit{hauteur} in \cite{bourbaki:comalg1-7} VI \S4.4, is the Krull dimension $\dim \Spec(R_w)$ of its valuation ring $R_w$. 

\smallskip

\begin{center}
\begin{tabular}[t]{|c|||c|c|c|c|c|} \hline
     type &  value group & $\bQ$-rank & rank & on $k$ & residue field \\ \hline \hline \hline
     0 & 1 & 0& 0 & trivial & $K$ \\ \hline \hline
     1h & $\bZ$ & 1& 1 & trivial & finite over $k$ \\ \hline
     1v & $\bZ$ & 1& 1 & $v$ & function field over $\kkp$ of \\ 
            & & & & & transcendence degree $1$ \\ \hline \hline
     2h & $\bZ \oplus \bZ$ lex. & 2 & 2 & $v$ & finite over $\kkp$ \\ \hline
     2v & $\bZ \oplus \bZ$ lex.  & 2 & 2 & $v$ & finite over $\kkp$ \\ \hline
     2u$_{\rm n}$  & $\bZ \oplus \bZ \gamma \subset \bR$ & 2 & 1 & $v$ & finite over $\kkp$ \\ \hline
     2u$_{\rm sm}$  & $\bZ$ & 1 & 1 & $v$ & infinite, algebraic over $\kkp$ \\ \hline
     2u$_{\rm alt}$  & $\bigcup_n \bruch{1}{e_n}\bZ$ & 1 & 1 & $v$ & algebraic over $\kkp$ \\ 
            &   with $\lim  e_n = \infty$ & & & & \\ \hline
  \end{tabular}
\end{center}

\medskip

The residue field for a valuation $w$ of type 2u$_{\rm smooth}$ has to be algebraic over $\kappa$ of infinite degree. Indeed, otherwise the extension $\fok \prec R_w$ had finite residue degree $f = [\kappa(w):\kappa]$ and finite index of value groups $e = (w(K): v(k))$, which implies that $K$ as a $k$ vector space has $\dim_k K = ef$, a contradiction. In particular, if we ultimately pick smooth centers $x_w$ and the residue field degree $[\kappa(x_w):\kappa]$ remains finite, then we actually deal with a valuation of type 2h.

\subsection{Valuations of the universal cover}

From now on we fix a geometric generic point $\bar{\eta}:\Spec(\Omega) \to X$ of $X$ as base point. Let $\tilde{K}$ be the function field of the associated pointed universal pro-\'etale cover $\tilde{X}$ of $X$, i.e, $\tilde{K} \subset \Omega$ is the maximal algebraic extension of $K$ which is unramified over $X$. We conclude that $\pi_1(X,\bar{\eta})$ equals $\Gal(\tilde{K}/K)$.

\subsubsection{The Riemann--Zariski space of the universal cover}
The prolongation $\Val_\fok(\tld K)$ of $\Val_\fok(K)$ to
$\tld K$ endowed with the patch topology is a projective limit 
\[
\Val_\fok(\tilde{K}) \xrightarrow{\sim} \varprojlim_{K'} \Val_\fok(K')
\]
of the spaces $\Val_\fok(K')$ equipped with the patch topology, where $K'$ ranges over all finite intermediate extensions $K'/K$ in $\tilde{K}/K$.  As above, one has a 
homeomorphism
\[
{\rm center} \, : \ \Val_\fok(\tilde{K}) \xrightarrow{\sim} \varprojlim_{X' \subset \cX'} \cX'_\cons
\]
and the subset 
\[
\Val_v (\tilde{K}) =  \{ \tilde{w} \in \Val_\fok(\tilde{K}) \ ; \ \tilde{w}|_k = v \} \subset  \Val(\tilde{K}) 
\]
is a closed subset in the patch topology described by the condition that $\tilde{w}(\pi)>0$ for a uniformizer $\pi$ of $\fok$. Thus $\Val_v(\tilde{K})$ is a compact, Hausdorff, pro-finite space which furthermore is canonically a pro-finite limit
\begin{equation} \label{eq:centervalvtilde}
{\rm center} \, : \ \Val_v (\tilde{K}) \xrightarrow{\sim} \varprojlim_{K',\cX'} \cX'_{\kkp,\cons}
\end{equation}
 of the pro-finite spaces $\cX'_{\kkp,\cons}$, where $\cX'_{\kkp,\cons}$ is the 
 reduced special fibre of $\cX'$ endowed with the constructible topology.

\subsubsection{Types and the universal cover}
The canonical restriction map $\Val_\fok(K') \to \Val_\fok(K)$ is surjective, and for $w' \mapsto w$, by the fundamental inequality, the residue field extension $\kappa(w')/\kappa(w)$ is finite and the inclusion of value groups $w(K) \subset w'(K')$ has finite index, see \cite{bourbaki:comalg1-7} VI \S8.
Hence the type of a valuation is preserved under the restriction map $\Val_\fok(K') \to \Val_\fok(K)$, and the classification into types also applies to valuations in $\Val_\fok(\tilde{K})$.

\subsubsection{Notational convenience}
The map $\Val_\fok(\tilde{K}) \to \Val_\fok(K)$ will be denoted by $\tilde{w} \mapsto w=\tilde{w}|_K$ which implicitly could also imply a choice of a preimage $\tilde{w}$ of the valuation $w$ if the latter happens to appear first.


\section{Unramified \textit{Hilbert Zerlegungstheorie}}  \label{sec:hilb}

We keep the notation and assumptions from Appendix~\ref{sec:val}.

\subsection{Nearby points} \label{sec:nearby}
For a geometric point $y$ on a model $\cX$ we set $\cX_y^{\rh} = \Spec (\OO^{\rh}_{\cX,y})$ for the \textbf{scheme of nearby points} and $\cX_y^\sh = \Spec (\OO^\sh_{\cX,y})$ for the \textbf{scheme of strictly nearby points}. The intersection with the generic fibre we denote by 
\[
\cU_y^{\rh} = \Spec (\OO^{\rh}_{\cX,y}\otimes_{\fok} k) \subset  \cX_y^{\rh} 
\qquad \text{ and } \qquad
 \cU_y^\sh = \Spec (\OO^\sh_{\cX,y}\otimes_{\fok} k) \subset  \cX_y^\sh.
 \]
For $y$ equal to the center $x_w$ of a valuation $w \in \Val_\fok(K)$, more precisely, for a choice of geometric point  above the closed point of the valuation ring which induces a geometric point $\ov{x}_w$ above each center,  we abbreviate 
\[
\cU_w^{\rh} := \cU_{\bar{x}_w}^{\rh} \subseteq \cX_w^{\rh} := \cX_{\bar{x}_w}^{\rh}
\qquad \text{ and } \qquad
\cU_w^\sh := \cU_{\bar{x}_w}^\sh \subseteq \cX_w^\sh := \cX_{\bar{x}_w}^\sh.
\]
In the limit over all models $\cX$ of $K$ we get 
\[
U_w^{\rh} = \varprojlim_\cX \cU_w^{\rh} \subseteq X_w^{\rh} = \varprojlim_\cX \cX_w^{\rh} 
\qquad \text{ and } \qquad
U_w^\sh = \varprojlim_\cX \cU_w^\sh \subseteq X_w^\sh = \varprojlim_\cX \cX_w^\sh.
\]
We note that $U_w^{\rh}$ (resp.\ $U_w^\sh$) is a limit of affine $p$-adic curves over $k$ (resp.\ $k^\nr$). In particular, the cohomological dimension of $U_w^\sh$ for \'etale constructible sheaves is at most $2$.

\subsection{Hilbert decomposition and inertia group}
Let us fix a choice of a geometric generic point $\bar{\xi}_y$ of $\cU_y^\sh$ such that $(\cU_y^\sh,\bar{\xi}_y) \to (X,\bar{\eta})$ becomes a pointed map.

The \textbf{decomposition group} {resp.\ \textbf{inertia group}} in the sense of Hilbert at $y$ is given by the image $\rD_y$, resp.\ $\rI_y$, of the natural map $\pi_1(\cU_y^{\rh},\bar{\xi}_y) \to \pi_1(X,\bar{\eta})$,  resp.\ $\pi_1(\cU_y^\sh,\bar{\xi}_y) \to \pi_1(X,\bar{\eta})$,  induced by the inclusions. We suppress the choice of base points in the notation for decomposition and inertia groups.

\subsection{Decomposition and inertia group of a valuation} 
For $w \in \Val_\fok(K)$ let  $\bar{\xi}_{w}$  be a geometric generic point of $U_w^\sh$ such that $(U_w^\sh,\bar{\xi}_w) \to (X,\bar{\eta})$ becomes a pointed map. 
The compatibility of $\bar{\xi}_w$ with $\bar{\eta}$ describes a unique prolongation $\tilde{w}$ of $w$ to $\tilde{K}$ by the property $K_{\tilde{w}}^\sh = \tilde{K} \cdot K_w^\sh$ and similarly $K_{\tilde{w}}^{\rh} = \tilde{K} \cdot K_w^{\rh}$  in $\Omega$. Here $K_w^{\rh}$ (resp.\ $K_w^\sh$) is a (strict) henselisation of $K$ in $w$, and similarly for $\tilde{w}$. We easily observe the following lemma.

\begin{lem}
For a valuation $w \in \Val_\fok(K)$ of \valht $2$ but not of type 2h we have $\Spec(K_w^{\rh}) = U_w^{\rh}$, whereas for $w$ type 2h refining $\alpha$ of type 1h the nearby points $U_w^{\rh}$ equals the spectrum of the valuation ring the exention of $\alpha$ to $K_w^{\rh}$ and moreover equals $U_\alpha^{\rh} = X_\alpha^{\rh}$. \hfill $\square$
\end{lem}

The \textbf{decomposition group} (resp.\ \textbf{inertia group}) in the sense of valuation theory of  $w$, or more precisely the prolongation $\tilde{w}|w$ to a valuation of $\tilde{K}$,
 is given by the image $\rD_{\tilde{w}|w}$  of 
$\pi_1(U_w^{\rh},\bar{\xi}_w) \to \pi_1(X,\bar{\eta})$, resp.\  the image $\rI_{\tilde{w}|w}$ of $\pi_1(U_w^\sh,\bar{\xi}_w) \to \pi_1(X,\bar{\eta})$.  

The dependence on  $\tilde{w}$ is through the choice of a path connecting the base points $\bar{\xi}_w$ and $\bar{\eta}$ to the effect of conjugating $\rD_{\tilde{w}|w}$ and $\rI_{\tilde{w}|w}$ within $\pi_1(X,\bar{\eta})$. If no confusion arises, we will simplify the notation to $\rD_w = \rD_{\tilde{w}|w}$ (resp.\  $\rI_w = \rI_{\tilde{w}|w}$).

\subsection{Reconciliation of valuation theory and arithmetic geometry}
The two viewpoints of inertia and decomposition groups are related via the compliance of $\pi_1$ with affine projective limits.  We may assume that  $\bar{\xi}_w$ induces $\bar{\xi}_{\bar{x}_w}  \in \cU_w^\sh$ for every model of $X$, and then find
\begin{equation} \label{eq:limitID}
\rD_{\tilde{w}|w} = \varprojlim_\cX \rD_{x_w} \quad \text{ and } \quad \rI_{\tilde{w}|w} = \varprojlim_\cX \rI_{x_w} ,
\end{equation}
where the limits are in fact simply intersections of closed subgroups in $\pi_1(X,\bar{\eta})$.

Moreover, let $\alpha$ be a valuation of \valht 1 and $y$ a geometric point localised in a closed point of the divisor $Y_\alpha$ associated to $\alpha$ on a suitable model $\cX$. Then we have the following diagram when the corresponding geometric points are compatibly chosen.
\[
\xymatrix@M+1ex@R-2.5ex{ & & \Spec(\Omega) \ar[dll]_{\bar{\xi}_\alpha} \ar[dl]^{\bar{\xi}_y} \ar[ddd]^{\bar{\eta}} \\
U^\sh_\alpha \ar@{.>}[r] \ar[d] & \cU^\sh_y \ar[d] & \\
U^{\rh}_\alpha  \ar@{.>}[dr] & \cU^{\rh}_y \ar[d] & \\
 & \cU_y^{\nis_\alpha} \ar[r] &  X
}
\]
The scheme $ \cU_y^{\nis_\alpha} $ is the generic fibre of $ \cX_y^{\nis_\alpha}$ which is the maximal strict \'etale neighbourhood in between $\cX_y^{\rh} \to \cX$ which is Nisnevich at $\alpha$, i.e., such that the point $\alpha$ splits in the image of $\bar{\xi}_\alpha$ after an appropriate choice is fixed.
With $\rD_{y,\alpha} = \im\Big(\pi_1(\cU_y^{\nis_\alpha},\bar{\xi}_y) \to \pi_1(X,\bar{\eta})\Big)$, we find
\begin{equation} \label{eq:ID}
\xymatrix{ \rI_\alpha \ar@{}[r]|{\displaystyle \subseteq} \ar@{}[d]|{\bigcap}  & \rI_y \ar@{}[r]|{\displaystyle \subseteq} & \rD_y \ar@{}[d]|{\bigcap} & \\
\rD_\alpha \ar@{}[rr]|{\displaystyle \subseteq} &  & \rD_{y,\alpha} \ar@{}[r]|{\displaystyle \subseteq} & \pi_1(X,\bar{\eta}) }
\end{equation}
Let $w \in \Val_\fok(K)$ be a valuation  of rational rank $2$, i.e. of type 2h or 2v, and a refinement of the valuation $\alpha$, then in the limit over all models we deduce from (\ref{eq:ID}) and (\ref{eq:limitID}) that
\[ \rI_\alpha \subseteq \rI_w \subseteq \rD_w \subseteq \rD_\alpha\]
because $\rD_\alpha = \varprojlim_{\cX} \rD_{x_w,\alpha}$. 

\end{appendix}


\end{document}